\documentclass{article}
\usepackage{amsfonts}

\usepackage{graphicx}
\usepackage{amsmath}

\numberwithin{equation}{section}
\newtheorem{theorem}{Theorem}

\newtheorem{corollary}[theorem]{Corollary}

\newtheorem{lemma}[theorem]{Lemma}

\newtheorem{proposition}[theorem]{Proposition}
\newtheorem{remark}[theorem]{Remark}

\newenvironment{proof}[1][Proof]{\textbf{#1.} }{\ \rule{0.5em}{0.5em}}
\input{tcilatex}

\begin{document}

\bigskip

\bigskip

\begin{center}
\textbf{ON THE EXISTENCE OF EXPONENTIALLY DECREASING SOLUTIONS OF THE
NONLINEAR LANDAU DAMPING PROBLEM.}

\bigskip

Hyung Ju Hwang\footnote{%
Department of Mathematics, Pohang University of Science and Technology,
Pohang 790-784, Republic of Korea.} and Juan J. L. Vel\'{a}zquez\footnote{%
ICMAT (CSIC-UAM-UC3M-UCM). Facultad de Matem\'{a}ticas, Universidad
Complutense, Madrid 28040, Spain.}

\bigskip
\end{center}

\textbf{Abstract.-} In this paper we prove the existence of a large class of
periodic solutions\ of the Vlasov-Poisson in one space dimension that decay
exponentially as $t\rightarrow \infty .$ The exponential decay is well known
for the linearized version of the Landau damping problem and it has been
proved in \cite{CM} for a class of solutions of the Vlasov-Poisson system
that behaves asymptotically as free streaming solutions and are sufficiently
flat in the space of velocities. The results in this paper enlarge the class
of possible asymptotic limits, replacing the flatness condition in \cite{CM}
by a stability condition for the linearized problem.

\bigskip


\bigskip

\section{Introduction.}

\bigskip

Landau damping is a remarkable property of collisionless plasmas. This
effect was discovered in \cite{Landau} and it consists in the exponential
damping of charge oscillations in the plasma due to the combined effect of
the electrical fields generated by the charges and the dispersion in the
particle velocities.

\bigskip

In more mathematical terms, this effect is usually studied using the
Vlasov-Poisson system for negatively charged particles moving in a constant
background of positive charges that makes the whole system electrically
neutral. A linearized Vlasov-Poisson problem is then derived assuming that
the inhomogeneities of the charge distribution are small. For rather general
initial distributions of particles and initial velocities the corresponding
charge disturbances of this linear system decay exponentially fast.

The mathematical theory of the Vlasov-Poisson equation, including global
existence results for general initial data in space dimension $N\leq3$ has
been established in several papers (cf. \cite{Glassey}, \cite{LP}, \cite
{Pfaffelmoser}).

Landau damping has been also extensively studied in the mathematical and
physical literature (cf. \cite{BohmGross}, \cite{Degond}, \cite{GuoStrauss}, 
\cite{Jackson}, \cite{Saenz}, \cite{Schaefer}). The results in these papers
hold only for the linearized version of the problem. The derivation of the
exponential decay for the linearized Vlasov-Poisson system, in the cases
where such a decay takes place, relies heavily on the analytic properties of
the initial data (cf. \cite{Landau}, \cite{Jackson}, \cite{Saenz}). However,
it has been proved in \cite{Schaefer} that exponential decay cannot be
expected, even for the linearized problem, for nonanalytic initial data.

There are also physical arguments that explain why the energy of the charge
waves tends to be transferred to the distribution of velocities of the
particles (cf. \cite{BohmGross}, \cite{Jackson}), yielding exponential decay
of the charge disturbances. In these arguments the linearity of the problem
plays an essential role, in particular it is possible to decompose an
arbitrary disturbance in monochromatic waves using Fourier transform.
However, it is not obvious how to generalize these physical arguments to the
whole nonlinear problem, even for small perturbations of the linear case.
Nevertheless, the numerical simulations in \cite{Guo} indicate that
exponential decay of the charge disturbances might be expected for the whole
nonlinear Vlasov-Poisson model and a large class of analytic initial data.

A rigorous mathematical proof of the existence of solutions to the
Vlasov-Poisson system in a circle for which the electrical field $E$
decreases exponentially has been obtained in \cite{CM}. However, a general
stability condition ensuring such exponential decay has not been obtained
yet.

We now precise the mathematical problem that we will be considered in this
paper. Assuming that there is a background of charge $n_{0}$ to make the
system electrically neutral, the mathematical problem that is usually
considered in the study of Landau damping is the following one: 
\begin{align}
f_{t}+vf_{x}+Ef_{v} & =0\;\;,\;\;-\infty<v<\infty \;\;,\;\;0<x<L\;\;,\;\;t>0
\label{S0E1} \\
f\left( 0,v,t\right) & =f\left( L,v,t\right)  \label{S0E2} \\
\phi_{xx}\left( x,t\right) & =\rho\left( x,t\right) =\int_{-\infty
}^{\infty}f\left( x,v,t\right) dv-n_{0}\;\;,\;\;0<x<L\;\;,\;\;t>0
\label{S0E3} \\
\phi\left( 0,t\right) & =\phi\left( L,t\right) \;\;,\;\;\phi_{x}\left(
0,t\right) =\phi_{x}\left( L,t\right) \;\;,\;\;t>0  \label{S0E4} \\
E & =\phi_{x}\;\;,\;\;0<x<L\;\;,\;\;t>0  \label{S0E5} \\
f\left( x,v,t\right) & =f_{0}\left( x,v\right) =f_{e}(v)+g_{0}\left(
x,v\right) ,\text{ \ \ }t=0   \label{S0E6}
\end{align}

where $f_{e}\left( v\right) $ satisfies: 
\begin{equation}
\int_{-\infty}^{\infty}f_{e}\left( v\right) dv=n_{0}   \label{E2}
\end{equation}

The physical basis on assuming the existence of the background of charge $%
n_{0}$ is that the system contains two types of particles with very
different masses. The lighter ones that can move easily are described by
means of the distribution $f\left( x,v,t\right) .$ On the other hand, there
are some heavier ''ions''\ that cannot move so easily and are replaced by
the charge background $n_{0}.$

Without loss of generality we can assume, using suitable units; 
\begin{equation}
n_{0}=1\;\;,\;\;L=2\pi   \label{E3}
\end{equation}

Notice that (\ref{E2}) implies that $f_{e}$ is a stationary solution of (\ref
{S0E1})-(\ref{S0E5}). Moreover, (\ref{S0E6}) and (\ref{E2}) imply: 
\begin{equation}
\int_{0}^{2\pi}\int_{-\infty}^{\infty}g_{0}dxdv=0   \label{E3a}
\end{equation}

In this paper we will obtain a large class of periodic solutions of the
nonlinear Vlasov-Poisson system (\ref{S0E1})-(\ref{S0E6}) in a bounded
domain for which the corresponding charge density $\rho $ and electric field 
$E$ decay exponentially as $t\rightarrow \infty .$ The result is a
perturbative one, in the sense that the derived solutions will be close to
solutions of the linearized equation. The results are close in spirit to the
ones in \cite{CM}. The main difference between the results in this paper and
the results in \cite{CM} is that we replace the condition on the limit
distribution in Lemma 3.1 by a stability condition on the limit distribution
that is equivalent to the stability conditions obtained in studies in the
physical literature on the Landau damping problem.

\bigskip

The main difficulty of the problem under consideration, even with the above
mentioned smallness conditions, is that a naive linearization near an
equilibrium distribution $f_{e}\left( v\right) ,$ say a maxwellian
distribution, does not allow to show that the remaining nonlinear terms are
small. This problem will be explained in detail in Section \ref{classical}.
The rationale behind this problem is that, although the electrical field and
the charge density converge to zero, the distribution function $f\left(
x,v,t\right) $ for the particles does not converge necessarily to the
equilibrium distribution $f_{e}\left( v\right) .$ On physical grounds this
might not be expected in the absence of collisions. Actually, the
dissipation of the energy contained in the field $E\left( x,t\right) $ must
result in the gain of kinetic energy of the particles of the system, as it
was noticed in \cite{BohmGross} (cf. also \cite{Jackson} for a clear
explanation of this). However, in the absence of dissipative mechanisms
there is no reason to expect for the long time distribution of particle
velocities to approximate to the equilibrium distribution $f_{e}\left(
v\right) $ or even to be spatially homogeneous. In fact, if the field $%
E\left( x,t\right) $ vanishes fast enough as $t\rightarrow\infty,$ the only
restriction that we can expect for the long time asymptotics of $f\left(
x,v,t\right) $ is to behave like a ''free streaming''\ function $%
f_{\infty}\left( x-vt,v\right) .$

The previous discussion suggests not to linearize around the equilibrium
distribution, but near a free streaming function $f_{\infty}\left(
x-vt,v\right) .$ It turns out that under suitable analyticity assumptions on 
$f_{\infty}$ and suitable smallness conditions ensuring that $f_{\infty}$ is
close to a stable equilibrium $f_{e}\left( v\right) ,$ it is possible to
obtain solutions of the Vlasov-Poisson system defined in $t\in\left(
0,\infty\right) ,$ such that $f\left( x,v,t\right) \sim f_{\infty}\left(
x-vt,v\right) $ as $t\rightarrow\infty,$ and $E\left( x,t\right) \rightarrow0
$ as $t\rightarrow\infty.$

It is interesting to notice that the linearization near the ''free
streaming''\ distribution $f_{\infty}\left( x-vt,v\right) $ has some
peculiar features. Indeed, the charge density $\rho_{\infty}\left(
x,t\right) \equiv\int f_{\infty}\left( x-vt,v\right) dv-1$ is not
identically zero for any $t\geq0$, as it could be expected from the fact
that $f\rightarrow f_{\infty}$ and $E\rightarrow0$ as $t\rightarrow\infty.$
It turns out that $\rho_{\infty}\left( x,t\right) $ approaches to zero
exponentially fast as $t\rightarrow\infty$ under suitable analyticity
assumptions on $f_{\infty}$ that will be made precise later. The decay to
zero of the density for the free streaming problem, due to the dispersion of
the velocities, is a well known fact (cf. \cite{Arnold}). The argument in
this paper indicates that this convergence to zero plays a crucial role in
the Landau damping problem. However, the averaging property is not the only
ground explaining Landau damping, since some additional stability conditions
must be satisfied by the distribution $f_{e}$. The conditions on $f_{e}$
that are required to obtain Landau damping have been explained by means of
intuitive physical arguments (cf. \cite{BohmGross}, \cite{Jackson}).
However, the role played by the averaging of the asymptotic distribution of
particle velocities has not been considered, to our knowledge, in the
physical study of this problem.

\bigskip

The main result of this paper is the following:

\bigskip

\begin{theorem}
\label{result}If $f_{e}\left( v\right) $ satisfies a suitable stability
condition for a linearized version of the Vlasov-Poisson problem there exist
initial data $f_{0}\left( x,v\right) =f_{e}\left( v\right) +g_{0}\left(
x,v\right) $ with $g_{0}$ satisfying (\ref{E3a}) such that the corresponding
solution of the system (\ref{S0E1})-(\ref{S0E6}) satisfies: 
\begin{equation*}
\left| E\left( x,t\right) \right| +\left| \rho \left( x,t\right) \right|
\leq Ce^{-\gamma t}\;\;,\;\;0\leq t<\infty 
\end{equation*}
for some $\gamma >0,\;C>0.$
\end{theorem}

\bigskip More precise conditions on $f_{e},\;g_{0}$ will be given in Section 
\ref{Nonlinear}, and Theorem \ref{Texist}.

The plan of the paper is the following. We recall in Section 2 some general
properties of the classical linearized Vlasov-Poisson system. In particular,
the difficulties that arise in trying to handle nonlinear terms with this
linearization procedure will be explained. Section 3 describes formally the
new linearization procedure. Section 4 formulates in detail a set of
conditions for $f_{e},\;g_{\infty}$ that will play a crucial role in the
rest of the paper. Section 5 analyses in detail mathematical properties of
the linearized operator. Section 6 contains the proof of exponential decay
for the nonlinear Vlasov-Poisson problem, in particular the proof of Theorem 
\ref{result}.

\bigskip

\section{\label{classical}The classical Landau linearized problem.}

In this Section we describe the linearized problem that was originally
considered by Vlasov (cf. \cite{Vlasov}) and Landau (cf. \cite{Landau}) and
that has been studied also in the mathematical and physical literature with
different levels of rigor \cite{BohmGross}, \cite{Degond}, \cite{Jackson}, 
\cite{Saenz}, \cite{Schaefer}). The goal is to explain shortly the kind of
difficulties that arise in trying to derive exponential decay for the field $%
E$ associated to the nonlinear Vlasov-Poisson problem taking as a starting
point the linearized problem. More precisely, we will explain why a naive
linearization argument does not allow to derive exponential decay for the
field. A linearization procedure that allows to obtain such exponential
decay estimates for the nonlinear problem will be introduced in Section \ref
{Lininfinity}. The procedure that we will use to solve the linearized
problem is equivalent to the one used by Landau, although it is different,
but closer to the approach that we will use in the forthcoming sections to
study the nonlinear problem.

We write: 
\begin{equation*}
f\left( x,v,t\right) =f_{e}\left( v\right) +g\left( x,v,t\right) 
\end{equation*}
We assume that $g$ is a small perturbation of the steady solution $%
f_{e}\left( v\right) .$ Neglecting formally quadratic terms in $g$ we obtain
the following problem that is the one usually considered in the study of the
linearized Landau damping problem:

\begin{align}
g_{t}+vg_{x}+Ef_{e,v} & =0\;\;,\;\;0<x<L\;\;,\;\;t>0\;,\;v\in\mathbb{R}
\label{L1} \\
\phi_{xx} & =\int gdv\;\;,\;\;E=\phi_{x}\;\;,\;\;0<x<L\;\;,\;\;t>0 
\label{L2} \\
g\left( x,v,0\right) & =g_{0}\left( x,v\right) \;\;,\;\;0<x<L\;\;,\;\;v\in%
\mathbb{R}  \label{L3} \\
\phi\left( 0,t\right) & =\phi\left( L,t\right) \;\;,\;\;\phi_{x}\left(
0,t\right) =\phi_{x}\left( L,t\right) \;\;,\;\;t>0  \label{L4} \\
g\left( 0,v,t\right) & =g\left( L,v,t\right) \;\;,\;v\in\mathbb{R\;}\text{\ }%
,~~t>0   \label{L5}
\end{align}

Then, using Duhamel's principle to solve (\ref{L1}), (\ref{L3}) we obtain: 
\begin{equation*}
g\left( x,v,t\right) =g_{0}\left( x-vt,v\right) -\int_{0}^{t}E\left(
x-v\left( t-s\right) ,s\right) f_{e,v}\left( v\right) ds 
\end{equation*}
where we assume that $g_{0},\;E$ are extended periodically to $x\in\mathbb{R}
$ with period $L.$

Plugging this identity into (\ref{L2}) we obtain the following
integro-differential equation for the electric field: 
\begin{align*}
E_{x}\left( x,t\right) & =\int_{-\infty}^{\infty}g\left( x,v,t\right)
dv=\int_{-\infty}^{\infty}g_{0}\left( x-vt,v\right) dv- \\
& -\int_{0}^{t}\left[ \int_{-\infty}^{\infty}E\left( x-v\left( t-s\right)
,s\right) f_{e,v}\left( v\right) dv\right] ds
\end{align*}

It is easier to study this equation using Fourier variables: 
\begin{align}
E\left( x,t\right) & =\sum_{n=-\infty}^{\infty}b_{n}\left( t\right) e^{inx}
\label{F1} \\
g_{0}\left( x,v\right) & =\sum_{n=-\infty}^{\infty}g_{n}\left( v\right)
e^{inx}  \notag
\end{align}
where the periodicity of $\phi$ requires $b_{0}\left( t\right) =0.$

Then the linearized problem becomes: 
\begin{equation*}
inb_{n}\left( t\right) =\int_{-\infty}^{\infty}dvg_{n}\left( v\right)
e^{-invt}-\int_{0}^{t}b_{n}\left( s\right) \left[ \int_{-\infty}^{\infty
}e^{-inv\left( t-s\right) }f_{e,v}\left( v\right) dv\right] ds 
\end{equation*}

We write: 
\begin{equation}
K\left( \xi\right) =\int_{-\infty}^{\infty}e^{-i\xi v}f_{e,v}\left( v\right)
dv   \label{L10}
\end{equation}

Then:

\begin{equation}
inb_{n}\left( t\right) +\int_{0}^{t}b_{n}\left( s\right) K\left( n\left(
t-s\right) \right) ds=G_{n}\left( t\right) \equiv\int_{-\infty}^{\infty
}dvg_{n}\left( v\right) e^{-invt}   \label{L7}
\end{equation}

In order to solve this equation we introduce a family of fundamental
solutions given by: 
\begin{equation}
inB_{n}\left( t;t_{0}\right) +\int_{0}^{t}B_{n}\left( s;t_{0}\right) K\left(
n\left( t-s\right) \right) ds=\delta\left( t-t_{0}\right)   \label{L8}
\end{equation}
where $t_{0}\geq0.$ It then follows that: 
\begin{equation}
b_{n}\left( t\right) =\int_{0}^{\infty}B_{n}\left( t;s\right) G_{n}\left(
s\right) ds   \label{L11}
\end{equation}

This family of convolution integral equations (\ref{L8}) can be solved using
Laplace transforms. Given a function $\varphi\left( t\right) $ in $\left\{
t>0\right\} $ we define: 
\begin{equation*}
\tilde{\varphi}\left( z\right) =\int_{0}^{\infty}\varphi\left( t\right)
e^{-zt}dt 
\end{equation*}

We recall that the Laplace transform of the convolution: 
\begin{equation*}
\left( \varphi\ast\psi\right) \left( t\right) =\int_{0}^{t}\varphi\left(
t-s\right) \psi\left( s\right) ds 
\end{equation*}
is given by: 
\begin{equation*}
\widetilde{\left( \varphi\ast\psi\right) }\left( t\right) =\tilde{\varphi }%
\left( z\right) \tilde{\psi}\left( z\right) 
\end{equation*}

On the other hand the Laplace transforms of the sequence of functions $%
K_{n}\left( t\right) =K\left( nt\right) $ are given by: 
\begin{equation*}
\widetilde{K_{n}}\left( t\right) =\int_{0}^{\infty}K\left( nt\right)
e^{-zt}dt=\frac{1}{n}\tilde{K}\left( \frac{z}{n}\right) 
\end{equation*}

Taking the Laplace transform of (\ref{L8}) it then follows that: 
\begin{equation*}
in\widetilde{B_{n}}\left( z;t_{0}\right) +\frac{1}{n}\tilde{K}\left( \frac{z%
}{n}\right) \widetilde{B_{n}}\left( z;t_{0}\right) =e^{-zt_{0}}
\end{equation*}

Then: 
\begin{equation}
\widetilde{B_{n}}\left( z;t_{0}\right) =\frac{e^{-zt_{0}}}{in+\frac{1}{n}%
\tilde{K}\left( \frac{z}{n}\right) }   \label{L9}
\end{equation}

Notice that the convolution structure of (\ref{L8}) implies that $%
B_{n}\left( t;s\right) =B_{n}\left( t-s;0\right) .$ This can be seen also
using (\ref{L9}). Moreover, classical properties of analytic functions imply
that the absence of zeros of $Q_{n}\left( z\right) \equiv in+\frac{1}{n}%
\tilde {K}\left( \frac{z}{n}\right) $ in the half-plane $\left\{ \func{Re}%
\left( z\right) >-\gamma\right\} $ for some $\gamma>0$ implies the
exponential decay of $B_{n}\left( t;0\right) $ as $t\rightarrow\infty$. The
function $Q_{n}\left( z\right) $ is, up to linear changes of variables, the
so-called Landau function, has been studied extensively in several papers.
Notice that an equivalent way of writing it, using (\ref{L10}), is the
following one:

\begin{align*}
\tilde{K}\left( z\right) & =\int_{0}^{\infty}d\xi e^{-\xi z}\int_{-\infty
}^{\infty}e^{-i\xi v}f_{e,v}\left( v\right) dv=\int_{-\infty}^{\infty}\frac{%
f_{e,v}\left( v\right) dv}{z+iv}=i\int_{-\infty}^{\infty}\frac {f_{e}\left(
v\right) dv}{\left( z+iv\right) ^{2}} \\
Q_{n}\left( z\right) & =i\left( n+\frac{1}{n}\int_{-\infty}^{\infty}\frac{%
f_{e}\left( v\right) dv}{\left( \frac{z}{n}+iv\right) ^{2}}\right)
\end{align*}

The following results can be derived from the results in the paper \cite
{Saenz} (cf. also \cite{Jackson} for related, although more formal results).

\begin{theorem}
\label{LinStab}Suppose that the function $f_{e}\left( v\right) $ is analytic
in the strip $\left| \func{Im}\left( v\right) \right| <A,\;g_{0}\left(
x,v\right) $ is analytic in $\left| \func{Im}\left( v\right) \right|
<A,\;\left| \func{Im}\left( x\right) \right| <A$ and satisfy 
\begin{equation*}
\left| f_{e}\left( v\right) \right| \leq\frac{B}{1+\left| v\right| ^{\alpha}}%
\;\;,\;\;\alpha>1\;\;,\;\;\left| \func{Im}\left( v\right) \right| <A 
\end{equation*}
\begin{equation*}
\left| g_{0}\left( x,v\right) \right| \leq\frac{B}{1+\left| v\right|
^{\alpha}}\;\;,\;\;\alpha>1\;\;,\;\;\left| \func{Im}\left( v\right) \right|
<A,\;\left| \func{Im}\left( x\right) \right| <A 
\end{equation*}
for some $A,B>0.$ Assume also that the function $Q_{n}\left( z\right) $ does
not have zeros in the half-plane $\func{Im}\left( z\right) \geq0$ for some
value of $n.$ Then the $n$-th Fourier coefficient of $E\left( x,t\right) $
defined in (\ref{F1}) decreases exponentially as $t\rightarrow \infty.$
\end{theorem}

\bigskip

\begin{theorem}
Suppose that $f_{e}\;$satisfy the assumptions in Theorem \ref{LinStab}.
Assume that the function $Q_{n}\left( z\right) $ has a zero in the
half-plane $\func{Im}\left( z\right) >0$ for some $n\in\mathbb{Z}.$ Then,
there exists $g_{0}\left( x,v\right) $ satisfying the assumptions in Theorem 
\ref{LinStab} such that: 
\begin{equation*}
\left\| E\left( \cdot,t\right) \right\| _{L^{p}\left( 0,2\pi\right) }\geq
Ce^{\gamma t}\;\;,\;\;1\leq p\leq\infty 
\end{equation*}
for some $\gamma>0,\;C>0.$
\end{theorem}

\bigskip

\begin{remark}
\label{StIn}Notice that the assumption on the zeros of $Q_{n}\left( z\right) 
$ in Theorem \ref{LinStab} is not satisfied for all the nonnegative initial
distributions $f_{e}\left( v\right) .$ There are several examples of such
functions yielding instabilities in \cite{Jackson}, \cite{Saenz}\ . For
instance, an example of function satisfying the analyticity assumptions on
Theorem \ref{LinStab} and where the corresponding function $Q_{n}\left(
z\right) $ has zeros in the half-plane $\func{Im}\left( z\right) >0$ for $%
n\neq0$ is (cf. \cite{Saenz}): 
\begin{equation*}
f_{e}\left( v\right) =\frac{4a^{5/2}}{3\pi^{\frac{1}{2}}}v^{4}e^{-av^{2}}
\end{equation*}

There exist distributions $f_{e}\left( v\right) $ for which the
corresponding functions $Q_{n}\left( z\right) $ yield stability. For
instance, the maxwellian distribution $f_{e}\left( v\right) =\frac
{e^{-v^{2}}}{\sqrt{\pi}}$ has been extensively studied in the physical
literature and it has been shown to satisfy the assumption in Theorem \ref
{LinStab} (cf. \cite{Jackson} , \cite{Saenz}). On the other hand there is a
proof in \cite{Jackson} of the fact that all the distributions of $%
f_{e}\left( v\right) $ with only one maximum yield stability for the
corresponding linearized problem.

An example, not so relevant physically, but where the condition on the roots
of $Q_{n}\left( z\right) $ in Theorem \ref{LinStab} can be checked easily
is: 
\begin{equation*}
f_{e}\left( v\right) =\frac{1}{1+v^{2}}
\end{equation*}

In this case, the function $Q_{n}\left( z\right) $ can be computed using
residues: 
\begin{equation*}
Q_{n}\left( z\right) =i\left[ n+\frac{1}{n}\frac{\pi}{\left( \frac{z}{n}%
+sgn\left( \func{Re}\left( \frac{z}{n}\right) \right) \right) ^{2}}\right] 
\end{equation*}
the zeroes of $Q_{n}\left( z\right) $ are: 
\begin{equation*}
z=-\left| n\right| \pm\frac{\sqrt{\pi}}{n}i 
\end{equation*}
\end{remark}

\bigskip

The following result can be obtained using the inversion formula for the
Laplace transform, as well as the methods in \cite{Saenz}.

\begin{theorem}
Suppose that $f_{e}$ satisfies the assumptions in Theorem \ref{LinStab} and
the function $Q_{n}\left( z\right) $ does not have zeroes in $\func{Im}%
\left( z\right) \geq0.$ Then the function $B_{n}\left( t;t_{0}\right) $
defined by means of (\ref{L8}) satisfies: 
\begin{equation*}
\left| B_{n}\left( t;t_{0}\right) \right| \leq\frac{C}{n}e^{-\gamma\left(
t-t_{0}\right) }\;\;,\;\;t\geq t_{0}\;\;,\;B_{n}\left( t;t_{0}\right)
=0\;\;,\;\;t<t_{0}
\end{equation*}
\end{theorem}

This estimate for the fundamental solution $B_{n}\left( t;t_{0}\right) $
implies an exponential decay for $E$ and for initial data $g_{0}$ satisfying
the assumptions in Theorem \ref{LinStab}. Indeed, these assumptions on $g_{0}
$ imply: 
\begin{equation*}
\left| G_{n}\left( t\right) \right| \leq\left| \int dvg_{n}\left( v\right)
e^{-invt}\right| \leq C_{n}\varepsilon e^{-\frac{A}{2}\left| n\right| t}
\end{equation*}
where the constants $C_{n}$ decrease faster than any power law if the
function $g_{0}$ is assumed to be $C^{\infty}$ in $x.$

Then, using (\ref{L11}): 
\begin{equation*}
\left| b_{n}\left( t\right) \right| \leq C_{n}\varepsilon e^{-bt}
\end{equation*}
for some $b>0,$ whence the exponential decay of $E$ follows.

\bigskip

Let us now indicate the difficulty that arises in trying to use the results
for the linearized problem in order to derive decay of the solutions for the
whole nonlinear problem. Suppose that we keep the neglected quadratic terms
in the linearized equation (\ref{L1}). Then: 
\begin{equation}
g_{t}+vg_{x}+Ef_{e,v}+Eg_{v}=0\;\;,\;\;0<x<L\;\;,\;\;t>0\;,\;v\in\mathbb{R} 
\label{L6}
\end{equation}

Using Duhamel's formula it then follows that: 
\begin{align*}
g\left( x,v,t\right) & =g_{0}\left( x-vt,v\right) -\int_{0}^{t}E\left(
x-v\left( t-s\right) ,s\right) f_{e,v}\left( v\right) ds \\
& -\int_{0}^{t}E\left( x-v\left( t-s\right) ,s\right) g\left( v\right) ds
\end{align*}

Plugging this identity into (\ref{L2}) as before it follows that: 
\begin{align*}
E_{x}\left( x,t\right) & =\int g\left( x,v,t\right) dv=\int g_{0}\left(
x-vt,v\right) dv- \\
& -\int_{0}^{t}\left[ \int E\left( x-v\left( t-s\right) ,s\right)
f_{e,v}\left( v\right) dv\right] ds- \\
& -\int_{0}^{t}\left[ \int E\left( x-v\left( t-s\right) ,s\right) g\left(
x,v,s\right) dv\right] ds
\end{align*}

Notice that, due to the presence of the last term containing the nonlinear
part of the integration $\int_{0}^{t}\left[ ...\right] ds$ it might not be
possible to ensure that this term converges to zero exponentially fast.
Actually, this does not happen in general. Therefore, it might not be
possible to treat this term as a small quadratic perturbation of the above
considered problem. This difficulty is not resolved by using more
sophisticated linearization procedures, for instance approximating the
characteristics of the whole nonlinear problem by the free streaming
characteristics plus a corrective term and expanding by using the Taylor
series in powers of the corrective term. The reason behind the failure of
these approximation procedures is that the function $g\left( x,v,t\right) $
does not converge to zero, although the electric field $E$ does. As a
consequence, all these linearization procedures contain terms that may not
be expected to be small as $t\rightarrow\infty$ and some method for handling
them must be found.

\section{\label{Lininfinity}Linearization near infinity: Formal computation.}

In this Section we describe at a formal level a more convenient way of
linearizing the Vlasov-Poisson system in order to derive estimates for the
nonlinear terms. The key idea is to linearize around the expected
asymptotics of the solutions as $t\rightarrow\infty.$ However, the function $%
f\left( x,v,t\right) $ is expected to be oscillatory as $t\rightarrow\infty,$
since in the absence of the field its dynamics would be given by the
transport equation, and then we might expect: 
\begin{equation*}
f\left( x,v,t\right) \sim F\left( x-vt,v\right) \;\;\text{as\ \ }%
t\rightarrow\infty 
\end{equation*}
for some function $F$. In order to obtain a function converging to a limit
as $t\rightarrow\infty$ it is convenient to introduce a new set of variables:

\begin{align}
z & =x-vt  \label{U1E1} \\
\tilde{f}\left( x-vt,v,t\right) & =f\left( x,v,t\right)   \label{U1E2}
\end{align}

Then: 
\begin{equation*}
f_{t}=\tilde{f}_{t}-v\tilde{f}_{z}\;\;,\;\;f_{x}=\tilde{f}_{z}\;\;,\;\;f_{v}=%
\tilde{f}_{v}-t\tilde{f}_{z}
\end{equation*}

Then (\ref{S0E1}) becomes: 
\begin{equation*}
\tilde{f}_{t}\left( z,v,t\right) -tE\left( z+vt,t\right) \tilde{f}_{z}\left(
z,v,t\right) +E\left( z+vt,t\right) \tilde{f}_{v}\left( z,v,t\right) =0 
\end{equation*}

On the other hand we can rewrite (\ref{S0E3}) as: 
\begin{equation*}
E_{x}\left( x,t\right) =\int_{-\infty}^{\infty}dvf\left( x,v,t\right)
-1=\int_{-\infty}^{\infty}\tilde{f}\left( x-wt,w,t\right) dw-1 
\end{equation*}
where we have used that $n_{0}=1.$ Also: 
\begin{equation*}
E_{z}\left( z,t\right) =\int_{-\infty}^{\infty}\tilde{f}\left(
z-wt,w,t\right) dw-1 
\end{equation*}

From now on, we will drop the tilde from $\tilde{f}$ in order to simplify
the notation. Therefore we need to study the problem:

\begin{equation}
f_{t}\left( z,v,t\right) -tE\left( z+vt,t\right) f_{z}\left( z,v,t\right)
+E\left( z+vt,t\right) f_{v}\left( z,v,t\right) =0   \label{N1}
\end{equation}
\begin{equation}
E_{z}\left( z,t\right) =\int_{-\infty}^{\infty}f\left( z-wt,w,t\right) dw-1 
\label{N2}
\end{equation}

More precisely, we will construct solutions $f\left( z,v,t\right) $
satisfying: 
\begin{equation}
f\left( z,v,t\right) \rightarrow f_{e}\left( v\right) +g_{\infty}\left(
z,v\right) \equiv f_{\infty}\left( z,v\right) \;\;\text{as\ \ }%
t\rightarrow\infty   \label{N3}
\end{equation}
where the functions $f_{e}\left( v\right) ,\;g_{\infty}\left( x,v\right) $
will be assumed to satisfy suitable analyticity assumptions that will be
made precise later, and the function $g_{\infty}\left( x,v\right) $ is
assumed to be of order of the small parameter $\varepsilon$ and: 
\begin{equation}
\int_{0}^{2\pi}g_{\infty}\left( z,v\right) dz=0   \label{N3a}
\end{equation}

The characteristic equations associated to (\ref{N1}) satisfy: 
\begin{equation}
Z\left( t,t;z,v\right) =z\;\;,\;\;V\left( t,t;z,v\right) =v   \label{N5a}
\end{equation}
and: 
\begin{equation}
\frac{\partial Z}{\partial s}\left( s,t;z,v\right) =-sE\left( Z+Vs,s\right)
\;\;,\;\;\frac{\partial V}{\partial s}\left( s,t;z,v\right) =E\left(
Z+Vs,s\right)   \label{N5}
\end{equation}

We define the functions: 
\begin{equation}
Z\left( \infty,t;z,v\right) =Z_{\infty}\left( t;z,v\right) \;\;,\;\;V\left(
\infty,t;z,v\right) =V_{\infty}\left( t;z,v\right)   \label{N5b}
\end{equation}
Then, the solution of (\ref{N1}) can be written as: 
\begin{equation*}
f\left( z,v,t\right) =f_{\infty}\left( Z_{\infty}\left( t;z,v\right)
,V_{\infty}\left( t;z,v\right) \right) 
\end{equation*}

Using (\ref{N3}) it follows that: 
\begin{equation*}
f\left( z,v,t\right) =f_{e}\left( V_{\infty}\left( t;z,v\right) \right)
+g_{\infty}\left( Z_{\infty}\left( t;z,v\right) ,V_{\infty}\left(
t;z,v\right) \right) 
\end{equation*}

Therefore: 
\begin{equation*}
E_{z}\left( z,t\right) =\int_{-\infty}^{\infty}\left[ f\left(
z-wt,w,t\right) -f_{e}\left( w\right) \right] dw 
\end{equation*}
whence: 
\begin{align}
E_{z}\left( z,t\right) & =\int_{-\infty}^{\infty}\left[ f_{e}\left(
V_{\infty}\left( t;z-wt,w\right) \right) -f_{e}\left( w\right) \right] dw
\label{N4} \\
& +\int_{-\infty}^{\infty}g_{\infty}\left( Z_{\infty}\left( t;z-wt,w\right)
,V_{\infty}\left( t;z-wt,w\right) \right) dw  \notag
\end{align}

In order to linearize (\ref{N4}) we argue as follows. We are trying to
construct solutions where $E$ decreases exponentially as $%
t\rightarrow\infty. $ Moreover, the field can be expected to be of order $%
\varepsilon.$ In particular the functions $Z,\;V$ can be expected to be
nearly constant. We can then approximate the solution of the equations (\ref
{N5}) to the leading order as: 
\begin{align}
Z_{\infty}\left( t;z,v\right) -z & =-\int_{t}^{\infty}sE\left( z+vs,s\right)
ds  \label{U1E3} \\
V_{\infty}\left( t;z,v\right) -v & =\int_{t}^{\infty}E\left( z+vs,s\right)
ds   \label{U1E4}
\end{align}

Linearizing the functions $f_{e},g_{\infty}$ using Taylor's theorem, and
taking into account (\ref{U1E3}), (\ref{U1E4}) we obtain: 
\begin{align}
E_{z}\left( z,t\right) & =\int_{t}^{\infty}ds\int_{-\infty}^{\infty
}f_{e,v}\left( w\right) E\left( z-wt+ws,s\right)
dw+\int_{-\infty}^{\infty}g_{\infty}\left( z-wt,w\right) dw  \notag \\
& -\int_{t}^{\infty}sds\int_{-\infty}^{\infty}\frac{\partial g_{\infty}}{%
\partial z}\left( z-wt,w\right) E\left( z-wt+ws,s\right) dw+  \notag \\
& +\int_{t}^{\infty}ds\int_{-\infty}^{\infty}\frac{\partial g_{\infty}}{%
\partial w}\left( z-wt,w\right) E\left( z-wt+ws,s\right) dw   \label{N6}
\end{align}

Equation (\ref{N6}) is the linearized problem as $t\rightarrow\infty.$ Our
goal is to show that this problem provides a good approximation for the
solutions of the original problem (\ref{S0E1})-(\ref{S0E6}) as $t\rightarrow
\infty.$

\bigskip

It is relevant to notice that (\ref{N6}) contains terms that can be expected
to be of order $\varepsilon^{2},$ namely the last two ones. On the other
hand, the first two terms on the right hand side of (\ref{N6}) can be
expected to be of order $\varepsilon.$ The reason for keeping the last two
terms is that due to the presence of the term $s$ in the third term of (\ref
{N6}) it is not ''a priori'' clear if this term can be neglected for times $%
t\approx\frac {1}{\varepsilon}$ or larger. The last term in (\ref{N6}) has
been kept to ensure that for the resulting problem $\int_{0}^{2\pi}E_{z}%
\left( z,t\right) dz=0.$ Indeed, assuming that this identity is satisfied
and taking into account that $\int_{0}^{2\pi}E\left( z,t\right) dz=0$ as
well as (\ref{N3a}) it follows, after integrating (\ref{N6}) in $z\in\left[
0,2\pi\right] :$%
\begin{align}
0 & =-\int_{0}^{2\pi}dz\int_{t}^{\infty}sds\int_{-\infty}^{\infty}\frac{%
\partial g_{\infty}}{\partial z}\left( z-wt,w\right) E\left(
z-wt+ws,s\right) dw+  \notag \\
& +\int_{0}^{2\pi}dz\int_{t}^{\infty}ds\int_{-\infty}^{\infty}\frac{\partial
g_{\infty}}{\partial w}\left( z-wt,w\right) E\left( z-wt+ws,s\right) dw 
\label{N7}
\end{align}

\bigskip In order to check this inequality we notice that: 
\begin{equation}
\frac{Dg_{\infty}}{Dw}\left( z-wt,w\right) =-t\frac{\partial g_{\infty}}{%
\partial z}\left( z-wt,w\right) +\frac{\partial g_{\infty}}{\partial w}%
\left( z-wt,w\right)   \label{N8}
\end{equation}

Therefore, the right hand side of (\ref{N7}) can be written as: 
\begin{equation}
\int_{t}^{\infty}ds\int_{0}^{2\pi}dz\int_{-\infty}^{\infty}dwE\left(
z-wt+ws,s\right) \left[ \left( t-s\right) \frac{\partial g_{\infty}}{%
\partial z}\left( z-wt,w\right) +\frac{Dg_{\infty}}{Dw}\left( z-wt,w\right) %
\right]   \label{N8a}
\end{equation}
and, integrating by parts we can transform (\ref{N8a}) into: 
\begin{align*}
& -\int_{t}^{\infty}ds\int_{0}^{2\pi}dz\int_{-\infty}^{\infty}dwg_{\infty
}\left( z-wt,w\right) \cdot \\
& \cdot\left[ \left( t-s\right) E_{z}\left( z-wt+ws,s\right) -\left(
t-s\right) E_{z}\left( z-wt+ws,s\right) \right] =0
\end{align*}
whence (\ref{N7}) holds.

To conclude this section, we remark that there is a way of writing (\ref{N6}%
) where it becomes apparent that the last two terms are really small
perturbations. Indeed, using (\ref{N8}) we can rewrite (\ref{N6}) as: 
\begin{align*}
E_{z}\left( z,t\right) & =\int_{t}^{\infty}ds\int_{-\infty}^{\infty
}f_{e,v}\left( w\right) E\left( z-wt+ws,s\right)
dw+\int_{-\infty}^{\infty}g_{\infty}\left( z-wt,w\right) dw \\
& +\int_{t}^{\infty}\frac{s}{t}ds\int_{-\infty}^{\infty}\left[ \frac{%
Dg_{\infty}}{Dw}\left( z-wt,w\right) \right] E\left( z-wt+ws,s\right) dw+ \\
& +\int_{t}^{\infty}ds\int_{-\infty}^{\infty}\left( 1-\frac{s}{t}\right) 
\frac{\partial g_{\infty}}{\partial w}\left( z-wt,w\right) E\left(
z-wt+ws,s\right) dw
\end{align*}
and, after integrating by parts in the third equation: 
\begin{align}
E_{z}\left( z,t\right) & =\int_{t}^{\infty}ds\int_{-\infty}^{\infty
}f_{e,v}\left( w\right) E\left( z-wt+ws,s\right)
dw+\int_{-\infty}^{\infty}g_{\infty}\left( z-wt,w\right) dw  \label{N9} \\
& +\int_{t}^{\infty}\left( t-s\right) \frac{s}{t}ds\int_{-\infty}^{\infty
}g_{\infty}\left( z-wt,w\right) E_{z}\left( z-wt+ws,s\right) dw+  \notag \\
& +\int_{t}^{\infty}ds\int_{-\infty}^{\infty}\left( 1-\frac{s}{t}\right) 
\frac{\partial g_{\infty}}{\partial w}\left( z-wt,w\right) E\left(
z-wt+ws,s\right) dw  \notag
\end{align}

Suppose now that we study (\ref{N9}) in a space of functions satisfying $%
\left| E\left( z,t\right) \right| +\left| E_{z}\left( z,t\right) \right|
\leq Me^{-\gamma t}$ for some suitable $\gamma>0,\;M>0.$ Then, the third and
fourth terms on the right side of (\ref{N9}) can be estimated as: 
\begin{equation*}
C\varepsilon Me^{-\gamma t}
\end{equation*}

It then follows that, for a suitable choice of $\gamma,$ the last two terms
in (\ref{N9}) can be expected to be small perturbative terms. This fact will
be made rigorous in Section \ref{Nonlinear}.

Therefore, we expect to be able to approximate (\ref{N9}) as: 
\begin{equation}
E_{z}\left( z,t\right)
=\int_{t}^{\infty}ds\int_{-\infty}^{\infty}f_{e,v}\left( w\right) E\left(
z-wt+ws,s\right) dw+\int_{-\infty}^{\infty }g_{\infty}\left( z-wt,w\right)
dw   \label{N10}
\end{equation}

Equation (\ref{N10}) is basically equivalent to the linearized Landau
damping problem studied in Section \ref{classical}, except for the fact that
the linearization has been made at $t=\infty.$ We will study first some
properties of the linearized problem (\ref{N10}) using Fourier analysis. We
will then study the whole nonlinear problem (\ref{S0E1})-(\ref{S0E6}) using
a suitable functional framework.

\section{Analyticity assumptions for $f_{e}$ and $g_{\infty}$ and some
consequences.}

We make precise now the assumptions made on the functions $f_{e},\;g_{\infty}
$.

\bigskip

\textbf{Assumptions\ (A):}

\begin{itemize}
\item  The function $f_{e}\left( v\right) $ is analytic in the strip $\left| 
\func{Im}\left( v\right) \right| \leq A$ and it satisfies in that set 
\begin{equation}
\left| f_{e}\left( v\right) \right| \leq\frac{B}{1+\left| v\right| ^{\alpha}}%
\;\;,\;\;\alpha>1\;\;,\;\;\alpha\neq2   \label{A1}
\end{equation}

\item  The function $g_{\infty}\left( x,v\right) $ is analytic in the sets $%
\left| \func{Im}\left( x\right) \right| \leq A,\;\left| \func{Im}\left(
v\right) \right| \leq A$ and it satisfies in this set: 
\begin{equation}
\left| g_{\infty}\left( x,v\right) \right| \leq\frac{\varepsilon }{1+\left|
v\right| ^{\alpha}}\;\;,\;\;\alpha>1\;\;,\;\;\alpha\neq2   \label{A2}
\end{equation}

\item  The function $g_{\infty}\left( x,v\right) $ is periodic in the $x$
variable with period $2\pi$ and it satisfies: 
\begin{equation}
\int_{0}^{2\pi}g_{\infty}\left( x,v\right) dx=0   \label{A3}
\end{equation}
\end{itemize}

The requirement $\alpha\neq2$ seems at a first glance a bit artificial. This
assumption has been made in order to avoid the onset of logarithmic terms
that would introduce nonessential technical difficulties. Notice that we can
always assume that $\alpha\neq2$ reducing the value of $\alpha$ a bit if
needed.

There is a consequence of (\ref{A2}) that will be used repeatedly in the
following.

\begin{lemma}
\label{expdecay}Suppose that (\ref{A2}), (\ref{A3}) are satisfied. Then the
function 
\begin{equation}
H\left( z,t\right) =\int_{-\infty}^{\infty}g_{\infty}\left( z-wt,w\right) dw 
\label{T2}
\end{equation}
satisfies: 
\begin{equation}
\left| H\left( z,t\right) \right| \leq C\varepsilon e^{-\gamma
t}\;\;,\;\;z\in\left[ 0,2\pi\right]   \label{T1}
\end{equation}
where $\gamma>0$ can be chosen arbitrarily close to $A,\;\gamma<A$ and $C>0$
depends only on $\alpha,\gamma,A.$
\end{lemma}

\begin{proof}
Due to the analyticity properties of the function $g_{\infty}$ we can
rewrite $H\left( z,t\right) $ for $\left| \func{Im}\left( z\right) \right|
<A+At$ using contour deformation as: 
\begin{equation}
H\left( z,t\right) =\int_{\Gamma\left( z,t\right) }g_{\infty}\left(
z-wt,w\right) dw   \label{T3}
\end{equation}
where $\Gamma\left( z,t\right) =\left\{ w\in\mathbb{C}:\func{Im}\left(
w\right) =a\left( z,t\right) \right\} $ with $\left| \func{Im}\left(
z\right) -ta\left( z,t\right) \right| <A,$\newline
$\left| a\left( z,t\right) \right| <A.$ It is readily seen that for any $z$
satisfying $\left| \func{Im}\left( z\right) \right| <A+At$ it is possible to
choose such $a\left( z,t\right) ,$ although in a nonunique way. For
instance, given $z$ in $\left| \func{Im}\left( z\right) -rt\right| =\beta<A$
with $\left| r\right| <A,\;\left| \beta\right| <A$ we can choose $a\left(
z,t\right) =r.$ It then follows from the representation formula (\ref{T3})
that $H\left( z,t\right) $ is analytic in $\left| \func{Im}\left( z\right)
\right| <A+At.$ Moreover, using (\ref{A2}) it follows that 
\begin{equation}
\left| H\left( z,t\right) \right| \leq C\varepsilon   \label{T4}
\end{equation}
in $\left| \func{Im}\left( z\right) \right| <A+At.$

The periodicity of $g_{\infty}$ in $z$ implies that $H\left( z,t\right) $ is
periodic in $z$ for any $t\geq0$ and (\ref{A3}) implies 
\begin{equation*}
\int_{0}^{2\pi}H\left( z,t\right) dz=0 
\end{equation*}

We can then write $H\left( z,t\right) $ using the following Fourier series: 
\begin{equation}
H\left( z,t\right) =\sum_{n\neq0}\frac{e^{inz}}{2\pi}\int_{0}^{2\pi
}e^{-in\xi}H\left( \xi,t\right) d\xi   \label{T5}
\end{equation}

Using the analyticity properties of $H\left( \cdot,t\right) $ we can rewrite
the Fourier coefficients in (\ref{T5}) as: 
\begin{equation*}
\int_{0}^{2\pi}e^{-in\xi}H\left( \xi,t\right) d\xi=\int_{\left[ 0,2\pi\right]
+iL}e^{-in\xi}H\left( \xi,t\right) d\xi 
\end{equation*}
where $L\in\left( -A-At,A+At\right) .$ Choosing $sign\left( L\right)
=-sign\left( n\right) $ and choosing $L=\gamma+\gamma t$ with $\gamma<A$
arbitrarily close to $A$ it then follows from (\ref{T4}) that: 
\begin{equation*}
\left| \int_{0}^{2\pi}e^{-in\xi}H\left( \xi,t\right) d\xi\right| \leq2\pi
C\varepsilon e^{-\gamma\left| n\right| \left( t+1\right) }
\end{equation*}
and plugging this estimate into (\ref{T5}) we obtain (\ref{T1}).
\end{proof}

\section{Construction of the fundamental solution of the linearized problem
at $t=\infty.$}

\bigskip

Equation (\ref{N10}) as well as the form of the nonlinear terms that have
been neglected in the derivation of it suggest to study the following
problem: 
\begin{equation}
E_{z}\left( z,t\right)
=\int_{t}^{\infty}ds\int_{-\infty}^{\infty}f_{e,v}\left( w\right) E\left(
z-wt+ws,s\right) dw+h\left( z,t\right)   \label{N10a}
\end{equation}
where $h\left( z,t\right) $ is a bounded function decreasing sufficiently
fast as $t\rightarrow\infty$ and satisfying 
\begin{equation}
\int_{0}^{2\pi}h\left( z,t\right) dz=0   \label{N10b}
\end{equation}

In order to study this problem we will construct a fundamental solution
associated with it. More precisely we will derive, using Fourier analysis,
an explicit formula for a solution of: 
\begin{align}
G_{z}\left( z,t\right) & =\int_{t}^{\infty}ds\int_{-\infty}^{\infty
}f_{e,v}\left( w\right) G\left( z-wt+ws,s\right) dw+\left[ \delta\left(
z\right) -\frac{1}{2\pi}\right] \delta\left( t\right) ,\;  \label{E1G1} \\
t & \in\mathbb{R,\;}z\in\mathbb{R}  \notag
\end{align}
satisfying: 
\begin{equation}
G\left( z,t\right) =0\;\;,\;\;t>0\;,\;\;z\in\mathbb{R}   \label{E1G2}
\end{equation}
\begin{equation}
G\left( z,t\right) =G\left( z+2\pi,t\right)   \label{E1G3}
\end{equation}

Using (\ref{N10b}) as well as the invariance of the homogeneous part of (\ref
{E1G1}) under spatial and time translations we obtain that a solution $%
E\left( z,t\right) $ of (\ref{N10a}) can be written as: 
\begin{align}
E\left( z,t\right) & =\int_{-\infty}^{\infty}ds\int_{0}^{2\pi}d\xi G\left(
z-\xi,t-s\right) h\left( \xi,s\right)  \label{E1G3a} \\
& =\int_{t}^{\infty}ds\int_{0}^{2\pi}d\xi G\left( z-\xi,t-s\right) h\left(
\xi,s\right)  \notag
\end{align}
since, due to (\ref{N10b}): 
\begin{equation*}
h\left( z,t\right) =\int_{-\infty}^{\infty}ds\int_{0}^{2\pi}d\xi\left[
\delta\left( z-\xi\right) -\frac{1}{2\pi}\right] \delta\left( t-s\right)
h\left( \xi,s\right) 
\end{equation*}

There is a function that plays a crucial role in the whole theory of Landau
damping and that appears in slightly different forms in different papers
devoted to this subject. This function, that is usually referred as the
Landau function takes the following form in our setting: 
\begin{equation}
\Phi\left( \eta;n\right) =\int_{\mathbb{R}}\frac{f_{e,v}\left( w\right) }{%
w-\eta}dw-n^{2}   \label{LandauFunction}
\end{equation}

This function is defined in $\eta\in\mathbb{C}\setminus\mathbb{R}%
\;,\;n=\pm1,\pm2,...$. If $f_{e}$ satisfies (\ref{A1}), the function $\Phi$
can be extended analytically to the domain $\left\{ \func{Im}\left(
\eta\right) >-A\right\} $ for any $n=\pm1,\pm2,...$. It is worth mentioning
that the function $\Phi\left( \eta;n\right) $ is discontinuous for $\eta \in%
\mathbb{R}.\;$Indeed, due to the Plemej-Sokolski formula (cf. \cite{Ahlfors})%
\texttt{:} 
\begin{equation*}
\Phi\left( \eta_{0}+i0\right) -\Phi\left( \eta_{0}-i0\right) =2\pi
if_{e,v}\left( \eta_{0}\right) 
\end{equation*}

Therefore the analytic continuation of the function $\Phi$ to the domain $%
\left\{ \func{Im}\left( \eta\right) >-A\right\} $ is not given by the
integral formula (\ref{LandauFunction}). This is a well known property of
the Landau function (cf. \cite{Jackson}).

\bigskip

The main result of this Section is summarized in the following
Theorem:\bigskip

\begin{theorem}
\label{Linear} Suppose that $f_{e}$ satisfies (\ref{A1}). Let $0<\delta<A.$
Define two functions $\psi_{\pm}\left( \eta\right) $ by means of: 
\begin{equation}
\psi_{\pm}\left( \eta\right) =\int_{C_{\pm\delta}}\frac{f_{e,v}\left(
w\right) }{\eta\pm w}dw,   \label{E1V1}
\end{equation}
where $C_{\pm\delta}=\mathbb{R}\pm i\delta.$

We define a function $Q\left( z,t\right) $ by means of: 
\begin{align}
Q_{z}\left( z,t\right) & =\frac{1}{\left( 2\pi\right) ^{2}}%
\sum_{n=-\infty\;,\;n\neq0}^{\infty}\frac{e^{inz}}{n}\int_{-\infty}^{\infty }%
\frac{e^{i\left| n\right| \eta t}\psi_{sign\left( n\right) }\left(
\eta\right) }{\left( 1-\frac{\psi_{sign\left( n\right) }\left( \eta\right) }{%
\left| n\right| n}\right) }d\eta  \label{E1V2} \\
Q\left( z,t\right) & =\int_{0}^{z}Q_{z}\left( \xi,t\right) d\xi+\frac
{1}{2\pi}\int_{0}^{2\pi}\xi Q_{z}\left( \xi,t\right) d\xi   \label{E1V3}
\end{align}
as well as a tempered distribution $G$ by means of: 
\begin{align}
G_{z}\left( z,t\right) & =Q_{z}\left( z,t\right) +\delta\left( t\right) 
\left[ \sum_{\ell=-\infty}^{\infty}\delta\left( z+2\pi\ell\right) -\frac{1}{%
2\pi}\right]  \label{E1V4} \\
G\left( z,t\right) & =\int_{0}^{z}G_{z}\left( \xi,t\right) d\xi+\frac
{1}{2\pi}\int_{0}^{2\pi}\xi G_{z}\left( \xi,t\right) d\xi   \label{E1V5}
\end{align}

Suppose that the Landau function defined in (\ref{LandauFunction}) does not
have zeroes in the half-plane $\left\{ \func{Im}\left( \eta\right)
\geq0\right\} .$ Then, the series defining the function $Q_{z}$ in (\ref
{E1V2}) converges for any $z\in\mathbb{R},$ $t\in\mathbb{R.}$ The
distributions $G,\;Q$ are supported in the set $t\leq0$ and they are
periodic in $z$ with period $2\pi.$ We have $Q\in C^{\infty}\left( \mathbb{%
R\times }\left( -\infty,0\right) \right) .$ Moreover, the following
estimates hold: 
\begin{align}
\left| Q\left( z,t\right) \right| +\left| Q_{z}\left( z,t\right) \right| &
\leq\frac{C}{\left( 1+\left| \frac{z}{t}\right| ^{\alpha }\right) }%
\;\;,\;\;0\leq t\leq1\;\;,\;\;-\pi\leq z\leq\pi  \label{E1V5a} \\
\left| Q\left( z,t\right) \right| +\left| Q_{z}\left( z,t\right) \right| &
\leq Ce^{-a\left| t\right| }\;\;,\;\;t\geq1\;\;,\;\;z\in \mathbb{R} 
\label{E1V5b}
\end{align}
where $\alpha$ is as in (\ref{A1}) and the constants $C,\;a$ depend on $%
B,\;\alpha.$

The distribution $G$ solves (\ref{E1G1})-(\ref{E1G3}).\texttt{\ }
\end{theorem}

In order to prove Theorem \ref{Linear} we begin by deriving the analyticity
properties of the functions $\psi_{\pm}$ as well as some basic estimates for
them.

\begin{lemma}
\label{Le1}Suppose that $f_{e}$ satisfies (\ref{A1}) and let $\delta$ be as
in Theorem \ref{Linear}. Then the functions $\psi_{\pm}\left( \eta\right) $
defined in (\ref{E1V1}) are analytic in $\left\{ \func{Im}\left( \eta\right)
>-\delta\right\} $ and they satisfy the following estimates: 
\begin{equation}
\left| \psi_{\pm}\left( \eta\right) \right| \leq\frac{C}{1+\left|
\eta\right| ^{\beta}}\;\;,\;\;\func{Im}\left( \eta\right) >-\frac{\delta}{2}%
\;\;,   \label{E1V6}
\end{equation}
\begin{equation}
\left| D_{\eta}^{k}\psi_{\pm}\left( \eta\right) \right| \leq\frac{C_{k}}{%
1+\left| \eta\right| ^{\beta}}\;\;,\;\;\func{Im}\left( \eta\right) >-\frac{%
\delta}{2}\;,\;\;k=1,2,...   \label{E1V8}
\end{equation}
where $1<\beta=\min\left\{ \alpha,2\right\} -\varepsilon_{0}$ with $%
\varepsilon_{0}>0$ that might be chosen arbitrarily small$,\;C$ depends only
on $\alpha,\;B,\;\delta,\;\varepsilon_{0}$ and $C_{k}$ depends on $%
\alpha,\;B,\;\delta,\;\varepsilon_{0},\;k$.
\end{lemma}

\begin{proof}
The analyticity of $\psi_{\pm}\left( \eta\right) $ in $\left\{ \func{Im}%
\left( \eta\right) >-\delta\right\} $ is just a consequence of the
analyticity of the functions $\frac{1}{\eta\pm w}$ for each $w\in
C_{\pm\delta}$ in the half-plane $\func{Im}\left( \eta\right) >-\delta$ as
well as the assumption (\ref{A1}).

To derive (\ref{E1V6}) we use the inequality: 
\begin{equation*}
\left| \frac{1}{\eta\pm w}-\frac{1}{\eta}\right| \leq\frac{\left| w\right| }{%
\left| \eta\right| \left| \eta\pm w\right| }
\end{equation*}

Then: 
\begin{equation}
\left| \psi_{\pm}\left( \eta\right) -\frac{1}{\eta}\int_{C_{\pm%
\delta}}f_{e,v}\left( w\right) dw\right| \leq\frac{B}{\left| \eta\right| }%
\int_{C_{\pm\delta}}\frac{\left| w\right| }{\left| \eta\pm w\right| }\frac{1%
}{1+\left| w\right| ^{\alpha}}\left| dw\right|   \label{E1V9}
\end{equation}

Using the fact that $\int_{C_{\pm\delta}}f_{e,v}\left( w\right) dw=0$ and
splitting the integral on the right side of (\ref{E1V9}) in the regions
where $\left| w\right| \leq\frac{\left| \eta\right| }{2}$ and $\left|
w\right| >\frac{\left| \eta\right| }{2}$ respectively, we obtain: 
\begin{align}
\left| \psi_{\pm}\left( \eta\right) \right| & \leq\frac{B}{\left|
\eta\right| ^{2}}\int_{C_{\pm\delta}\cap\left\{ \left| w\right| \leq \frac{%
\left| \eta\right| }{2}\right\} }\frac{\left| w\right| }{1+\left| w\right|
^{\alpha}}\left| dw\right| +  \label{E1V10} \\
& +\frac{B}{\left| \eta\right| }\int_{C_{\pm\delta}\cap\left\{ \left|
w\right| >\frac{\left| \eta\right| }{2}\right\} }\frac{1}{\left| \eta\pm
w\right| \left| w\right| ^{\alpha-1}}\left| dw\right|  \notag
\end{align}

We can estimate the first integral on the right side of (\ref{E1V10}) by a
constant if $\alpha>2$ and as $C\left| \eta\right| ^{-2+\alpha}$ if $%
\alpha<2.$ (Notice that $\alpha\neq2,$ and therefore the logarithmic case
does not occur). On the other hand we can estimate the second integral on
the right hand side of (\ref{E1V10}) introducing the rescaling $w=\left|
\eta\right| \zeta.$ Then: 
\begin{equation*}
\left| \psi_{\pm}\left( \eta\right) \right| \leq\frac{C}{\left| \eta\right|
^{\beta}}+\frac{B}{\left| \eta\right| ^{\alpha}}\int _{\frac{C_{\pm\delta}}{%
\left| \eta\right| }\cap\left\{ \left| \zeta\right| >\frac{1}{2}\right\} }%
\frac{1}{\left| \frac{\eta}{\left| \eta\right| }\pm\zeta\right| \left|
\zeta\right| ^{\alpha-1}}\left| d\zeta\right| 
\end{equation*}

The integral in this formula is convergent, since $\alpha>1.$ The main
contribution to this integral for large $\left| \eta\right| $ is due to the
region where $\left| \frac{\eta}{\left| \eta\right| }\pm\zeta\right| $ is
small. The closest distance between $\zeta$ and $\mp\frac{\eta}{\left|
\eta\right| }$ is of order $\frac{1}{\left| \eta\right| }.$ Therefore: 
\begin{equation*}
\int_{\frac{C_{\pm\delta}}{\left| \eta\right| }\cap\left\{ \left|
\zeta\right| >\frac{1}{2}\right\} }\frac{1}{\left| \frac{\eta}{\left|
\eta\right| }\pm\zeta\right| \left| \zeta\right| ^{\alpha-1}}\left|
d\zeta\right| \leq C_{\delta}\left[ \left| \log\left( \left| \eta\right|
\right) \right| +1\right] 
\end{equation*}
whence, choosing $\varepsilon_{0}>0$ arbitrarily small, (\ref{E1V6}) follows.

Estimate (\ref{E1V8}) can be proved combining (\ref{E1V6}) as well as the
fact that the functions $\psi_{\pm}\left( \eta\right) $ are analytic in $%
\left\{ \func{Im}\left( \eta\right) >-\delta\right\} $ and the classical
Cauchy's inequalities for analytic functions (cf. \cite{Ahlfors}).
\end{proof}

\bigskip

The functions $\psi_{\pm}\left( \eta\right) $ are closely related to the
Landau function $\Phi$ (cf. (\ref{LandauFunction})). We reformulate some
properties of the functions $\psi_{\pm}\left( \eta\right) $ in terms of
properties of the function $\Phi$ that have been studied often in the
literature (cf. for instance \cite{BohmGross}).

\bigskip

\begin{lemma}
\label{Le4}Suppose that $f_{e}$ satisfies (\ref{A1}). Suppose that the
Landau function $\Phi\left( \eta;n\right) $ defined in (\ref{LandauFunction}%
) does not have zeroes in the region $\left\{ \func{Im}\left( \eta\right)
\geq0\right\} $ for any $n=\pm1,\pm2,...$. Then, there exist $\nu_{0}>0$ and 
$\theta>0$ such that 
\begin{equation}
\left| 1-\frac{\psi_{sign\left( n\right) }\left( \eta\right) }{\left|
n\right| n}\right| \geq\theta   \label{J1a}
\end{equation}
for $\func{Im}\left( \eta\right) >-\nu_{0}.$
\end{lemma}

\begin{proof}
Using (\ref{E1V1}) we have: 
\begin{equation}
1-\frac{\psi_{sign\left( n\right) }\left( \eta\right) }{\left| n\right| n}=1-%
\frac{1}{n^{2}}\int_{C_{\pm\delta}}\frac{f_{e,v}\left( w\right) }{%
w+sign\left( n\right) \eta}dw\equiv-\frac{1}{n^{2}}\Phi\left( -sign\left(
n\right) \eta;n\right)   \label{B1}
\end{equation}

Notice that (\ref{LandauFunction}) as well as the fact that $f_{e}\left(
w\right) $ takes real values for $w\in\mathbb{R}$ implies: 
\begin{equation}
\Phi\left( \bar{\eta};n\right) =\overline{\left( \Phi\left( \eta;n\right)
\right) }   \label{B2}
\end{equation}

Therefore, if $\Phi\left( \eta;n\right) $ does not have zeroes for $\left\{ 
\func{Im}\left( \eta\right) \geq0\right\} ,$ it does not have zeroes for $%
\left\{ \func{Im}\left( \eta\right) \leq0\right\} $ either. Moreover, due to
(\ref{B2}) we can restrict our attention to the case $sign\left( n\right)
=-1.$ Due to\ (\ref{A1}) we can extend analytically the function $\Phi\left(
\eta;n\right) $ to the domain $\left\{ \func{Im}\left( \eta\right)
>-A\right\} $ by means of the formula: 
\begin{equation*}
\Phi\left( \eta;n\right) =\int_{\mathbb{R}-Ai}\frac{f_{e,v}\left( w\right) }{%
w-\eta}dw-n^{2}
\end{equation*}

It then follows, due to (\ref{A1}) that for $\left| \eta\right| >\rho$ with $%
\rho$ independent of $n$: 
\begin{equation*}
\left| \frac{1}{n^{2}}\Phi\left( -sign\left( n\right) \eta;n\right) \right|
\geq\frac{1}{2}
\end{equation*}

On the other hand, since $\Phi\left( \eta;n\right) $ does not have zeroes
for $\eta\in\left\{ \func{Im}\left( \eta\right) \geq0\right\} $, $\left|
\eta\right| \leq\rho$ it follows by continuity that the analytic extension
of $\frac{1}{n^{2}}\Phi\left( \eta;n\right) $ to $\left\{ \func{Im}\left(
\eta\right) >-A\right\} $ does not have zeroes in the region $-\nu_{0}\leq%
\func{Im}\left( \eta\right) \leq0,\;\left| \eta\right| \leq\rho$ for some $%
\nu_{0}>0$ sufficiently small. Notice that $\nu_{0}$ can be chosen uniformly
in $n.$ Using (\ref{B1}) the result follows.
\end{proof}

\begin{remark}
The absence of zeros of the Landau function $\Phi\left( \eta;n\right) $ in
the half plane $\left\{ \func{Im}\left( \eta\right) \geq0\right\} $ for any $%
n=\pm1,\pm2,...$ is precisely the condition required for the stability of
the solutions of the linearized Landau problem studied in Section \ref
{classical}. Therefore, the functions $f_{e}\left( v\right) $ yielding
stability for the problem considered there yield also stability for the
linearized problem considered in this Section. In particular, the examples
of stability and instability in Remark \ref{StIn} are also valid for the
linearized problem considered in this Section.
\end{remark}

\bigskip

We now derive some estimates on the Fourier coefficients in the series (\ref
{E1V2}) that will ensure the convergence of the series.

\bigskip

From now on we will write by shortness $\pm$ instead of $sign\left( n\right)
.$

\bigskip

\begin{lemma}
\label{Le2}Suppose that $f_{e}$ satisfies (\ref{A1}) and let $\delta$ be as
in Theorem \ref{Linear}.\texttt{\ }Suppose that the Landau function $%
\Phi\left( \eta;n\right) $ defined in (\ref{LandauFunction}) does not have
zeroes in the region $\left\{ \func{Im}\left( \eta\right) \geq0\right\} $
for any $n=\pm1,\pm2,...$. Then: 
\begin{equation}
\left| \int_{-\infty}^{\infty}\frac{e^{i\left| n\right| \eta t}\psi_{\pm
}\left( \eta\right) }{\left( 1-\frac{\psi_{\pm}\left( \eta\right) }{\left|
n\right| n}\right) }d\eta\right| \leq\frac{C_{\gamma}}{\left| n\right|
^{\gamma}}\frac{1}{\left| t\right| ^{\gamma}}\;\;,\;t\in \mathbb{R\;},%
\mathbb{\;}\;\gamma>0\;,\;n\neq0   \label{Est1}
\end{equation}
where $C_{\gamma}$ depends on $\alpha,\;B,\;\delta,\;\gamma.$ Moreover, we
have: 
\begin{equation}
\int_{-\infty}^{\infty}\frac{e^{i\left| n\right| \eta t}\psi_{\pm}\left(
\eta\right) }{\left( 1-\frac{\psi_{\pm}\left( \eta\right) }{\left| n\right| n%
}\right) }d\eta=0\;\;\text{for\ \ }t>0   \label{Est2}
\end{equation}
and 
\begin{equation}
\left| \int_{-\infty}^{\infty}\frac{\psi_{\pm}\left( \eta\right) }{\left( 1-%
\frac{\psi_{\pm}\left( \eta\right) }{\left| n\right| n}\right) }d\eta\right|
\leq\frac{C}{n^{2}}\;\;,\;n\neq0   \label{Est3}
\end{equation}
where $C$ depends on $\alpha,\;B,\;\delta,\;\gamma.$
\end{lemma}

\begin{proof}
Notice that, for any $\ell=1,2,...$ we have: 
\begin{equation*}
e^{i\left| n\right| \eta t}=\frac{1}{\left( i\left| n\right| t\right) ^{\ell}%
}\frac{\partial^{\ell}}{\partial\eta^{\ell}}\left( e^{i\left| n\right| \eta
t}\right) 
\end{equation*}

Then, integrating by parts we obtain: 
\begin{equation}
\int_{-\infty}^{\infty}\frac{e^{i\left| n\right| \eta t}\psi_{\pm}\left(
\eta\right) }{\left( 1-\frac{\psi_{\pm}\left( \eta\right) }{\left| n\right| n%
}\right) }d\eta=\frac{\left( -1\right) ^{\ell}}{\left( i\left| n\right|
t\right) ^{\ell}}\int_{-\infty}^{\infty}e^{i\left| n\right| \eta t}\frac{%
\partial^{\ell}}{\partial\eta^{\ell}}\left( \frac{\psi_{\pm}\left(
\eta\right) }{\left( 1-\frac{\psi_{\pm}\left( \eta\right) }{\left| n\right| n%
}\right) }\right) d\eta   \label{J1}
\end{equation}

Using (\ref{J1a}) in Lemma \ref{Le4}, it then follows from (\ref{E1V8})
that: 
\begin{equation}
\left| \frac{\partial^{\ell}}{\partial\eta^{\ell}}\left( \frac{\psi_{\pm
}\left( \eta\right) }{\left( 1-\frac{\psi_{\pm}\left( \eta\right) }{\left|
n\right| n}\right) }\right) \right| \leq\frac{C_{\ell}}{1+\left| \eta\right|
^{\beta}}   \label{J2}
\end{equation}
with $\beta$ as in Lemma \ref{Le1}. Combining (\ref{J1}), (\ref{J1a}), (\ref
{J2}) we obtain 
\begin{equation}
\left| \int_{-\infty}^{\infty}\frac{e^{i\left| n\right| \eta t}\psi_{\pm
}\left( \eta\right) }{\left( 1-\frac{\psi_{\pm}\left( \eta\right) }{\left|
n\right| n}\right) }d\eta\right| \leq\frac{C_{\ell}}{\left( \left| n\right|
t\right) ^{\ell}}\;\;,\;\;t\in\mathbb{R}\;,\mathbb{\;}n\neq0\;,\;%
\ell=0,1,2,...   \label{J3}
\end{equation}

Notice that for $\ell=0$ (\ref{J3}) follows immediately from (\ref{E1V6}).
Estimate (\ref{Est1}) follows for $\gamma\in\left( \ell-1,\ell\right) $ by
interpolation.

The identity (\ref{Est2}) can be obtained using the fact that in the half
plane $\func{Im}\left( \eta\right) >0$ where the function $\frac{%
\psi_{\pm}\left( \eta\right) }{\left( 1-\frac{\psi_{\pm}\left( \eta\right) }{%
\left| n\right| n}\right) }$ is analytic and it is bounded as $\frac{C_{\ell}%
}{1+\left| \eta\right| ^{\beta}}$ we have also the estimate $\left|
e^{i\left| n\right| \eta t}\right| =e^{-\left| n\right| \func{Im}\left(
\eta\right) t}\leq1.$ Then (\ref{Est2}) follows by deforming the contour of
integration $\left( -\infty,\infty\right) $ to $\left( -\infty,\infty\right)
+iR$ with $R>0$ and taking the limit $R\rightarrow\infty.$

Finally, we prove (\ref{Est3}) using the identity: 
\begin{equation}
\int_{-\infty}^{\infty}\frac{\psi_{\pm}\left( \eta\right) }{\left( 1-\frac{%
\psi_{\pm}\left( \eta\right) }{\left| n\right| n}\right) }%
d\eta=\int_{-\infty}^{\infty}\psi_{\pm}\left( \eta\right) d\eta+\frac
{1}{\left| n\right| n}\int_{-\infty}^{\infty}\frac{\left( \psi_{\pm}\left(
\eta\right) \right) ^{2}}{\left( 1-\frac{\psi_{\pm}\left( \eta\right) }{%
\left| n\right| n}\right) }d\eta   \label{J4}
\end{equation}

The first term on the right hand side of (\ref{J4}) can be computed using (%
\ref{E1V1}): 
\begin{equation*}
\int_{-\infty}^{\infty}\psi_{\pm}\left( \eta\right) d\eta=\lim
_{R\rightarrow\infty}\int_{-R}^{R}\int_{C_{\pm\delta}}\frac{f_{e,v}\left(
w\right) }{\eta\pm w}dwd\eta=\lim_{R\rightarrow\infty}\int_{C_{\pm%
\delta}}dwf_{e,v}\left( w\right) \int_{-R}^{R}\frac{d\eta}{\eta\pm w}
\end{equation*}

Since $\int_{-R}^{R}\frac{d\eta}{\eta\pm w}=\log\left( \frac{R\pm w}{-R\pm w}%
\right) $ is bounded for large $R$ and $w\in C_{\pm\delta}$ we can apply
Lebesgue dominated convergence Theorem to obtain: 
\begin{equation}
\int_{-\infty}^{\infty}\psi_{\pm}\left( \eta\right) d\eta=\int_{C_{\pm
\delta}}dwf_{e,v}\left( w\right) \lim_{R\rightarrow\infty}\left[ \log\left( 
\frac{R\pm w}{-R\pm w}\right) \right] =-\pi i\int_{C_{\pm\delta
}}dwf_{e,v}\left( w\right) =0   \label{J5}
\end{equation}

Using (\ref{E1V6}), (\ref{J1a}) and (\ref{J5}) in (\ref{J4}) we obtain (\ref
{Est3}).
\end{proof}

\bigskip

\bigskip We can now prove the convergence of the series defining $%
Q_{z}\left( z,t\right) $ in (\ref{E1V2}).

\begin{lemma}
Suppose that the assumptions of Theorem \ref{Linear} are satisfied. The
series on the right hand side of (\ref{E1V2}) is convergent for any $t\in%
\mathbb{R}.$ The function $Q_{z}$ defined by means of (\ref{E1V2}) is
identically zero for $t>0,$ $Q_{z}\left( \cdot,t\right) \in C^{\infty}\left( 
\mathbb{R}\right) $ for any $t<0$ and $Q_{z}\left( \cdot,t\right) \in
C^{1}\left( \mathbb{R}\right) $ for $t=0.$ $Q_{z}\left( \cdot,t\right) $
defined by means of (\ref{E1V2}) is periodic with period $2\pi$ for any $t\in%
\mathbb{R}$ and it satisfies $\int_{0}^{2\pi}Q_{z}\left( z,t\right) dz=0.$
\end{lemma}

\begin{proof}
This Lemma is just a consequence of (\ref{Est1})-(\ref{Est3}).
\end{proof}

\bigskip

Notice that (\ref{Est1})-(\ref{Est3}) are not strong enough to derive
uniform estimates for the function $Q_{z}\left( z,t\right) $ or its
derivatives if $t$ is close to zero. This is made in the following Lemma,
where the self-similar structure of $Q_{z}\left( z,t\right) $ is derived.

\bigskip

\begin{lemma}
Suppose that the assumptions of Theorem \ref{Linear} are satisfied. Then 
\begin{equation*}
Q_{z}\left( z,t\right) =f_{e}\left( \frac{z}{t}\right) +R\left( z,t\right)
\;\;\;,\;\;z\in\mathbb{R}\;,\;\left| z\right| \leq \pi\;,\;\;-1\leq t<0 
\end{equation*}

where: 
\begin{equation*}
\left| R\left( z,t\right) \right| +\left| R_{z}\left( z,t\right) \right|
\leq C\;\;,\;\;z\in\mathbb{R}\;,\;\;t\in\mathbb{R}
\end{equation*}
\end{lemma}

\begin{proof}
We can write 
\begin{equation}
\int_{-\infty}^{\infty}\frac{e^{i\left| n\right| \eta t}\psi_{\pm}\left(
\eta\right) }{\left( 1-\frac{\psi_{\pm}\left( \eta\right) }{\left| n\right| n%
}\right) }d\eta=\int_{-\infty}^{\infty}e^{i\left| n\right| \eta
t}\psi_{\pm}\left( \eta\right) d\eta+\Omega_{n}\left( t\right)   \label{Cm1}
\end{equation}
where: 
\begin{equation}
\Omega_{n}\left( t\right) =\frac{1}{\left| n\right| n}\int_{-\infty
}^{\infty}\frac{e^{i\left| n\right| \eta t}\left( \psi_{\pm}\left(
\eta\right) \right) ^{2}}{\left( 1-\frac{\psi_{\pm}\left( \eta\right) }{%
\left| n\right| n}\right) }d\eta   \label{C0}
\end{equation}
Using (\ref{E1V6}), (\ref{J1a}), (\ref{Est3}) we have:

\begin{equation}
\left| \Omega_{n}\left( t\right) \right| \leq\frac{C}{n^{2}}   \label{C1}
\end{equation}

The term $\int_{-\infty}^{\infty}e^{i\left| n\right| \eta t}\psi_{\pm
}\left( \eta\right) d\eta$ can be computed explicitly. To this end we derive
a more convenient formula for this integral. Using (\ref{E1V1}) we obtain: 
\begin{align*}
& \int_{-\infty}^{\infty}e^{i\left| n\right| \eta t}\psi_{\pm}\left(
\eta\right) d\eta \\
& =\lim_{R\rightarrow\infty}\int_{-R}^{R}e^{i\left| n\right| \eta
t}\psi_{\pm}\left( \eta\right) d\eta \\
& =\lim_{R\rightarrow\infty}\int_{-R}^{R}e^{i\left| n\right| \eta
t}\int_{C_{\pm\delta}}\frac{f_{e,v}\left( w\right) }{\eta\pm w}dwd\eta \\
& =\lim_{R\rightarrow\infty}\int_{C_{\pm\delta}}dwf_{e,v}\left( w\right)
\int_{-R}^{R}\frac{e^{i\left| n\right| \eta t}}{\eta\pm w}d\eta \\
& =\int_{C_{\pm\delta}}dwf_{e,v}\left( w\right) \lim_{R\rightarrow\infty
}\int_{-R}^{R}\frac{e^{i\left| n\right| \eta t}}{\eta\pm w}d\eta
\end{align*}

We just need to compute the integral for $t<0.$ This computation can be made
by using residues. More precisely, we can replace the contour of integration 
$\left[ -R,R\right] $ by $\left[ -R,R\right] \cup\left\{ z:z=Re^{i\theta
},\;\theta\in\left[ -\pi,0\right] \right\} .$ The contribution to the
integral of the half-circle disappears as $R\rightarrow\infty,$ Indeed, the
integrand can be estimated as $\frac{C}{R}$ if $\left| \func{Im}\left(
\eta\right) \right| \leq A$ is of order one, and as $\frac {Ce^{-A\left|
nt\right| }}{R}$ for $\func{Im}\left( \eta\right) \leq-A.$ Using the fact
that the length of the half circle is bounded by $CR,$ it follows that the
contribution of the integral to the half-circle disappears taking the limits 
$R\rightarrow\infty,$ $A\rightarrow\infty.$ Thus:
\begin{equation*}
\lim_{R\rightarrow\infty}\int_{-R}^{R}\frac{e^{i\left| n\right| \eta t}}{%
\eta\pm w}d\eta=-2\pi ie^{\mp i\left| n\right| wt}=-2\pi ie^{-inwt}
\end{equation*}
whence, using also the analyticity properties of $f_{e}$: 
\begin{equation}
\int_{-\infty}^{\infty}e^{i\left| n\right| \eta t}\psi_{\pm}\left(
\eta\right) d\eta=-2\pi i\int_{-\infty}^{\infty}dwf_{e,v}\left( w\right)
e^{-inwt}   \label{C2}
\end{equation}

Using (\ref{Cm1}), (\ref{C2}) we can rewrite (\ref{E1V2}) as: 
\begin{align}
Q_{z}\left( z,t\right) & =\omega_{0}\left( z,t\right) +\omega_{1}\left(
z,t\right) \;\;,\;\;t<0  \label{C3} \\
\omega_{0}\left( z,t\right) & =-\frac{i}{\left( 2\pi\right) }%
\sum_{n=-\infty\;,\;n\neq0}^{\infty}\frac{e^{inz}}{n}\int_{-\infty}^{\infty
}dwf_{e,v}\left( w\right) e^{-inwt}  \notag \\
\omega_{1}\left( z,t\right) & =\frac{1}{\left( 2\pi\right) ^{2}}%
\sum_{n=-\infty\;,\;n\neq0}^{\infty}\frac{\Omega_{n}\left( t\right) }{n}%
e^{inz}  \notag
\end{align}

Notice that due to (\ref{C1}) the function $\omega_{1}\left( z,t\right) $
satisfies: 
\begin{equation}
\left| \omega_{1}\left( z,t\right) \right| +\left| \omega_{1,z}\left(
z,t\right) \right| \leq C\;\;,\;\;z\in\mathbb{R}\;,\;\;t\in\mathbb{R} 
\label{C4}
\end{equation}

On the other hand the function $\omega_{0}\left( z,t\right) $ can be
explicitly computed. Indeed, let us define the distribution: 
\begin{equation*}
T\left( z,t\right) \equiv-\frac{i}{\left( 2\pi\right) }\sum_{n=-\infty
\;,\;n\neq0}^{\infty}\frac{e^{in\left( z-wt\right) }}{n}
\end{equation*}

We have, in the sense of distributions: 
\begin{equation*}
\frac{\partial T}{\partial z}\left( z,t;w\right) =\frac{1}{2\pi}%
\sum_{n=-\infty\;,\;n\neq0}^{\infty}e^{inz}e^{-inwt}=\sum_{\ell=-\infty
}^{\infty}\delta\left( z-wt+2\pi\ell\right) -\frac{1}{2\pi}
\end{equation*}
whence, using the formula: 
\begin{equation*}
T\left( z,t;w\right) =\int_{0}^{z}T_{z}\left( \xi,t;w\right) d\xi+\frac
{1}{2\pi}\int_{0}^{2\pi}\xi T_{z}\left( \xi,t;w\right) d\xi 
\end{equation*}
we obtain: 
\begin{align*}
T\left( z,t;w\right) & =\sum_{\ell=-\infty}^{\infty}\left[ \chi\left(
z-wt+2\pi\ell\right) -\chi\left( -wt+2\pi\ell\right) \right] -\left( \frac{z%
}{2\pi}+\frac{1}{2}\right) + \\
& +\frac{1}{2\pi}\sum_{\ell=-\infty}^{\infty}\left( wt-2\pi\ell\right) \left[
\chi\left( 2\pi\left( \ell+1\right) -wt\right) -\chi\left(
2\pi\ell-wt\right) \right]
\end{align*}
where $\chi\left( \cdot\right) $ denotes the characteristic function whose
support is $\left( 0,\infty\right) .$

Writing 
\begin{equation*}
\omega_{0}\left( z,t\right) =\int_{-\infty}^{\infty}dwf_{e,v}\left( w\right)
T\left( z,t;w\right) 
\end{equation*}
it then follows that: 
\begin{equation*}
\omega_{0}\left( z,t\right) =\sum_{\ell=-\infty}^{\infty}\int_{\frac
{2\pi\ell}{t}}^{\frac{z+2\pi\ell}{t}}f_{e,v}\left( w\right) dw+\frac{1}{2\pi}%
\sum_{\ell=-\infty}^{\infty}\int_{\frac{2\pi\ell}{t}}^{\frac{2\pi\left(
\ell+1\right) }{t}}f_{e,v}\left( w\right) \left( wt-2\pi\ell\right) dw 
\end{equation*}
whence, after some integration: 
\begin{align*}
\omega_{0}\left( z,t\right) & =\sum_{\ell=-\infty}^{\infty}\left[
f_{e}\left( \frac{z+2\pi\ell}{t}\right) -f_{e}\left( \frac{2\pi\ell}{t}%
\right) \right] + \\
& +\frac{z}{2\pi}\sum_{\ell=-\infty}^{\infty}f_{e}\left( \frac{2\pi\left(
\ell+1\right) }{t}\right) -\frac{t}{2\pi}\sum_{\ell=-\infty}^{\infty}\int_{%
\frac{2\pi\ell}{t}}^{\frac{2\pi\left( \ell+1\right) }{t}}f_{e}\left(
w\right) dw
\end{align*}

Using (\ref{C3}) we obtain: 
\begin{equation*}
Q_{z}\left( z,t\right) =f_{e}\left( \frac{z}{t}\right) +R\left( z,t\right) 
\end{equation*}
\begin{align*}
R\left( z,t\right) & =-f_{e}\left( 0\right) +\sum_{\ell=-\infty,\ell
\neq0}^{\infty}\left[ f_{e}\left( \frac{z+2\pi\ell}{t}\right) -f_{e}\left( 
\frac{2\pi\ell}{t}\right) \right] + \\
& +\frac{z}{2\pi}\sum_{\ell=-\infty}^{\infty}f_{e}\left( \frac{2\pi\left(
\ell+1\right) }{t}\right) -\frac{t}{2\pi}\sum_{\ell=-\infty}^{\infty}\int_{%
\frac{2\pi\ell}{t}}^{\frac{2\pi\left( \ell+1\right) }{t}}f_{e}\left(
w\right) dw+\omega_{1}\left( z,t\right)
\end{align*}

Using (\ref{A1}) and (\ref{C4}) it follows that: 
\begin{equation*}
\left| R\left( z,t\right) \right| +\left| R_{z}\left( z,t\right) \right|
\leq C\;\;,\;\;z\in\mathbb{R}\;,\;\left| z\right| \leq \pi\;,\;-1\leq t<0 
\end{equation*}
Thus the result follows.
\end{proof}

We can now derive estimates for $Q_{z}\left( z,t\right) $ for $t$ large.

\begin{lemma}
Suppose that the assumptions of Theorem \ref{Linear} hold. Then: 
\begin{equation*}
\left| Q\left( z,t\right) \right| +\left| Q_{z}\left( z,t\right) \right|
\leq Ce^{at}
\end{equation*}
for $t<-1,$ with $a>0$ and $C>0$.
\end{lemma}

\begin{proof}
Combining (\ref{E1V6}) and (\ref{J1a}) we can write using contour
deformation: 
\begin{equation*}
\int_{-\infty}^{\infty}\frac{e^{i\left| n\right| \eta t}\psi_{\pm}\left(
\eta\right) }{\left( 1-\frac{\psi_{\pm}\left( \eta\right) }{\left| n\right| n%
}\right) }d\eta=\int_{\mathbb{R}-i\gamma}\frac{e^{i\left| n\right| \eta
t}\psi_{\pm}\left( \eta\right) }{\left( 1-\frac{\psi_{\pm }\left(
\eta\right) }{\left| n\right| n}\right) }d\eta 
\end{equation*}
thus we have 
\begin{equation*}
\left| \int_{-\infty}^{\infty}\frac{e^{i\left| n\right| \eta t}\psi_{\pm
}\left( \eta\right) }{\left( 1-\frac{\psi_{\pm}\left( \eta\right) }{\left|
n\right| n}\right) }d\eta\right| \leq Ce^{\gamma\left| n\right|
t}\;\;,\;\;t<0 
\end{equation*}
for some $\gamma>0.$ It then follows from (\ref{E1V2}) that $Q_{z}\left(
z,t\right) $ is periodic in $z$ with period $2\pi,$ as well as analytic and
bounded in $\left| \func{Im}\left( z\right) \right| \leq \frac{\gamma t}{2},$
$t<-1$. Arguing as in the proof of Lemma \ref{expdecay} we obtain that $%
Q_{z}\left( z,t\right) $ satisfies: 
\begin{equation*}
\left| Q_{z}\left( z,t\right) \right| \leq Ce^{at}\;\;,\;t\leq-1 
\end{equation*}
for some $a>0.$ A similar estimate for $Q\left( z,t\right) $ then follows
from the formula 
\begin{equation*}
Q\left( z,t\right) =\int_{0}^{z}Q_{z}\left( \xi,t\right) d\xi+\frac
{1}{2\pi}\int_{0}^{2\pi}\xi Q_{z}\left( \xi,t\right) d\xi 
\end{equation*}
and the result follows.
\end{proof}

To conclude the proof of Theorem \ref{Linear} it only remains to show that
the distribution $G$ defined in (\ref{E1V4}), (\ref{E1V5}) solves (\ref{E1G1}%
)-(\ref{E1G3}) in the sense of distributions. To this end we will
reformulate (\ref{E1G1}) in term of the Fourier transform of $G$ and we will
verify that the solution of the resulting equation is the one given by the
Fourier transform of $G$ that can be computed using (\ref{E1V2})-(\ref{E1V5}%
).

\bigskip

\begin{lemma}
\label{Ldist}Suppose that the assumptions of Theorem \ref{Linear} are
satisfied. Then the distribution $G\left( z,t\right) $ defined by means of (%
\ref{E1V4}), (\ref{E1V5}) solves (\ref{E1G1})-(\ref{E1G3}) in the sense of
distributions.
\end{lemma}

\begin{proof}
We expand $G\left( z,t\right) $ using Fourier series: 
\begin{equation}
G\left( z,t\right) =\sum_{n=-\infty}^{n=\infty}g_{n}\left( t\right) e^{inz} 
\label{S1F1}
\end{equation}

Plugging (\ref{S1F1}) into (\ref{E1G1}) we obtain, after some computations
that the functions $g_{n}\left( t\right) $ satisfy: 
\begin{equation}
ing_{n}\left( t\right) =\int_{t}^{\infty}dsg_{n}\left( s\right)
\int_{-\infty}^{\infty}f_{e,v}\left( w\right) e^{-inw\left( t-s\right) }dw+%
\frac{1}{2\pi}\delta\left( t\right) \;\;,\;\;n\neq0   \label{E1G4}
\end{equation}

In order to solve these equations we compute the Fourier transform of the
functions $g_{n}\left( t\right) $ that we define by means of: 
\begin{equation*}
\tilde{g}_{n}\left( \theta\right) =\frac{1}{\sqrt{2\pi}}\int_{-\infty
}^{\infty}g_{n}\left( t\right) e^{-it\theta}dt 
\end{equation*}

We define also the functions: 
\begin{equation}
\Phi_{\pm}\left( \theta\right) =\int_{-\infty}^{\infty}\varphi_{\pm}\left(
t\right) e^{-it\theta}dt   \label{S1F3}
\end{equation}
where: 
\begin{align*}
\varphi_{\pm}\left( t\right) & =\int_{-\infty}^{\infty}f_{e,v}\left(
w\right) e^{\mp iwt}dw\;\;,\;\;t<0 \\
\varphi_{\pm}\left( t\right) & =0\;\;,\;\;t\geq0
\end{align*}

The functions $\varphi_{\pm}$ decrease exponentially as $t\rightarrow-\infty.
$ This can be seen using (\ref{A1}) and contour deformation to derive the
following representation formula for these functions: 
\begin{equation}
\varphi_{\pm}\left( t\right) =\int_{C_{\pm\delta}}f_{e,v}\left( w\right)
e^{-iwt}dw\;\;,\;\;t<0   \label{E1F5}
\end{equation}
where $C_{\pm\delta}=\mathbb{R}\pm i\delta,$ $0<\delta<A.$

Applying the Fourier transform to (\ref{E1G4}) and using the fact that a
convolution is transformed to a product by this transformation we obtain: 
\begin{equation}
in\tilde{g}_{n}\left( \theta\right) =\frac{1}{\left| n\right| }\Phi_{\pm
}\left( \frac{\theta}{\left| n\right| }\right) \tilde{g}_{n}\left(
\theta\right) +\frac{1}{2\pi}\frac{1}{\sqrt{2\pi}}   \label{S1F2}
\end{equation}

On the other hand, using (\ref{E1F5}) we can rewrite the functions $%
\Phi_{\pm }$ in the following form: 
\begin{align*}
\Phi_{\pm}\left( \theta\right) & =\int_{-\infty}^{0}e^{-it\theta}\left[
\int_{C_{\pm\delta}}f_{e,v}\left( w\right) e^{\mp itw}dw\right]
dt=\int_{C_{\pm\delta}}f_{e,v}\left( w\right) \left[ \int_{-%
\infty}^{0}e^{-it\left( \theta\pm w\right) }dt\right] dw \\
& =i\int_{C_{\pm\delta}}\frac{f_{e,v}\left( w\right) }{\theta\pm w}dw
\end{align*}

Then (\ref{S1F2}) implies: 
\begin{equation*}
\tilde{g}_{n}\left( \theta\right) =\frac{1}{\sqrt{\left( 2\pi\right) ^{3}}}%
\frac{1}{i}\frac{1}{n-\frac{1}{\left| n\right| }\left[ \int_{C_{\pm \delta}}%
\frac{f_{e,v}\left( w\right) }{\frac{\theta}{\left| n\right| }\pm w}dw\right]
}
\end{equation*}
or, using the inversion formula for the Fourier transform: 
\begin{equation*}
g_{n}\left( t\right) =\frac{1}{\left( 2\pi\right) ^{2}i}\int_{-\infty
}^{\infty}\frac{e^{it\theta}d\theta}{n-\frac{1}{\left| n\right| }\left[
\int_{C_{\pm\delta}}\frac{f_{e,v}\left( w\right) }{\frac{\theta}{n}\pm w}dw%
\right] }\;\;,\;\;n\neq0 
\end{equation*}

It is readily seen that this formula is the same that the one for the
Fourier coefficients of the function $G$ in (\ref{E1V5}). Notice that the
choice (\ref{E1V5}) implies $g_{0}\left( t\right) =0.$ This concludes the
proof of the result.
\end{proof}

\bigskip

\bigskip

\begin{proof}[End of the Proof of Theorem \ref{Linear}]
It is just a consequence of Lemmas \ref{Le1}-\ref{Ldist}.
\end{proof}

Using Theorem \ref{Linear} we can prove the following result that will play
a crucial role in the analysis of the nonlinear problem:

\begin{proposition}
\label{regestimate}Suppose that the assumptions of Theorem \ref{Linear}
hold. Let $h\in C\left( \mathbb{R}\times\mathbb{R}^{+}\right) $ be a
function satisfying (\ref{N10b}), $h\left( z+2\pi,t\right) =h\left(
z,t\right) $ for $z\in\mathbb{R},\;t\geq0,$ as well as the estimate: 
\begin{equation*}
\left| h\left( z,t\right) \right| \leq Be^{-\gamma t}\;\;,\;\;z\in \mathbb{R}%
\;\;,\;t\geq0 
\end{equation*}
with $0<\gamma<a.$

Then there exists a function $E\left( z,t\right) \in C\left( \mathbb{R}\times%
\mathbb{R}^{+}\right) $ solving (\ref{N10a}), satisfying $\int
_{0}^{2\pi}E\left( z,t\right) dz=0,\;E\left( z+2\pi,t\right) =E\left(
z,t\right) $ and: 
\begin{equation*}
\left| E\left( z,t\right) \right| +\left| E_{z}\left( z,t\right) \right|
\leq CBe^{-\gamma t}
\end{equation*}
for some $C>0$ which depends only on $A,\;B,\;\alpha.$
\end{proposition}

\begin{proof}
The desired solution $E\left( z,t\right) $ can be obtained by means of the
formula (\ref{E1G3a}). Notice that we just need the values of $h\left(
z,t\right) $ in $t\geq0$ to obtain $E\left( z,t\right) $ for $t\geq0.$
Moreover, due to the linearity of the estimate we can assume $B=1.$ It is
readily seen that $\int_{0}^{2\pi}E\left( z,t\right) dz=0.$ Differentiating (%
\ref{E1G3a}) with respect to $z$ we obtain: 
\begin{equation*}
E_{z}\left( z,t\right) =\int_{-\infty}^{\infty}ds\int_{0}^{2\pi}d\xi
G_{z}\left( z-\xi,t-s\right) h\left( \xi,s\right) 
\end{equation*}

Using (\ref{E1V4}) it then follows that: 
\begin{equation*}
E_{z}\left( z,t\right) =h\left( z,t\right)
+\int_{-\infty}^{\infty}ds\int_{0}^{2\pi}d\xi Q_{z}\left( z-\xi,t-s\right)
h\left( \xi,s\right) 
\end{equation*}

Using (\ref{E1V5a}), (\ref{E1V5b}) it then follows that: 
\begin{equation*}
\left| E_{z}\left( z,t\right) \right| \leq Ce^{-\gamma t}
\end{equation*}

Finally we use the representation formula: 
\begin{equation*}
E\left( z,t\right) =\int_{0}^{z}E_{z}\left( \xi,t\right) d\xi+\frac
{1}{2\pi}\int_{0}^{2\pi}\xi E_{z}\left( \xi,t\right) d\xi 
\end{equation*}
to obtain: 
\begin{equation*}
\left| E\left( z,t\right) \right| \leq Ce^{-\gamma t}
\end{equation*}
whence the result follows.
\end{proof}

\section{\label{Nonlinear}Nonlinear problem: On the existence of
exponentially decaying solutions.}

In this Section we will prove the existence of a large class of solutions of
the nonlinear problem (\ref{S0E1})-(\ref{S0E6}) that decrease exponentially
fast as $t\rightarrow\infty.$

\begin{theorem}
\label{Texist}Suppose that $f_{e},\;g_{\infty}$ satisfy Assumptions (A).
Assume that the Landau function defined in (\ref{LandauFunction}) does not
have zeroes in the half-plane $\left\{ \func{Im}\left( \eta\right)
\geq0\right\} .$ Then there exist $\varepsilon_{0}=\varepsilon_{0}\left(
\alpha,A,B\right) >0,\;L=L\left( \alpha,A,B\right) $ such that for any $%
f_{e},\;g_{\infty}$ satisfying Assumptions (A) with $\varepsilon
\leq\varepsilon_{0}$ there exists a solution $f\left( z,v,t\right) $ of (\ref
{N1})-(\ref{N2}) defined for $z\in\mathbb{R},\;v\in\mathbb{R}$ and $0\leq
t<\infty$ and satisfying (\ref{N3}) as well as: 
\begin{equation*}
\left| E\left( x,t\right) \right| \leq C\varepsilon e^{-Lt}\;\;\text{for\ \ }%
0\leq t<\infty 
\end{equation*}
\end{theorem}

Using the change of variables (\ref{U1E1}), (\ref{U1E2}) we can prove the
following result:

\begin{corollary}
There exists a function $f_{0}\in C^{1}\left( \mathbb{R}\times\mathbb{R}%
\right) ,$ with $f_{0}\left( x+2\pi,v\right) =f_{0}\left( x,v\right) $ such
that the corresponding solutions of the system (\ref{S0E1})-(\ref{S0E6}) are
defined for $0\leq t<\infty$ and they satisfy: 
\begin{equation*}
\left| E\left( x,t\right) \right| \leq C\varepsilon e^{-Lt}\;\;\text{for\ \ }%
0\leq t<\infty 
\end{equation*}
\end{corollary}

\begin{remark}
In all the following $C$ will denote a numerical constant depending only on $%
\alpha,\;B,\;A$ that might change from line to line.
\end{remark}

\bigskip

\begin{proof}[Proof of Theorem \ref{Texist}]
The strategy is to solve (\ref{N4}) by means of a fixed point argument in a
space of functions satisfying: 
\begin{equation*}
\left\| E\right\| =\sup_{t>0}\left\{ e^{\gamma t}\left[ \left| E\left(
z,t\right) \right| +\left| E_{z}\left( z,t\right) \right| \right] \right\}
<\infty 
\end{equation*}

Suppose that $\left\| E\right\| <\infty.$ Then, integrating the
characteristic equations (\ref{N5}) it follows that: 
\begin{align*}
\left| V_{\infty}\left( t;z,w\right) -w\right| & \leq C\left\| E\right\|
e^{-\gamma t}\;\;,\;\;0\leq t<\infty\; \\
\left| Z_{\infty}\left( t;z,w\right) -z\right| & \leq C\left\| E\right\|
\left( t+1\right) e^{-\gamma t}\;\;,\;\;0\leq t<\infty
\end{align*}

Using Taylor's expansion in (\ref{N4}) we obtain: 
\begin{align*}
E_{z}\left( z,t\right) & =\int_{-\infty}^{\infty}f_{e,v}\left( w\right)
\left( V_{\infty}\left( t;z-wt,w\right) -w\right) dw+\int_{-\infty
}^{\infty}g_{\infty}\left( z-wt,w\right) dw+ \\
& +\int_{-\infty}^{\infty}\frac{\partial g_{\infty}}{\partial z}\left(
z-wt,w\right) \left( Z_{\infty}\left( t;z-wt,w\right) -\left( z-wt\right)
\right) dw+ \\
& +\int_{-\infty}^{\infty}\frac{\partial g_{\infty}}{\partial w}\left(
z-wt,w\right) \left( V_{\infty}\left( t;z-wt,w\right) -w\right) dw+R\left(
z,t\right)
\end{align*}
where: 
\begin{equation*}
\left| R\left( z,t\right) \right| \leq C\left\| E\right\| ^{2}\left(
t+1\right) ^{2}e^{-2\gamma t}
\end{equation*}

In order to derive better approximations to the differences $%
Z_{\infty}\left( t;z-wt,w\right) -\left( z-wt\right) ,\;V_{\infty}\left(
t;z-wt,w\right) -w$ we need to derive estimates for the solutions of the
characteristic equations: 
\begin{align}
\frac{\partial Z}{\partial s} & =-sE\left( Z+Vs,s\right)  \label{W2} \\
\frac{\partial V}{\partial s} & =E\left( Z+Vs,s\right)   \label{W3}
\end{align}
\begin{equation}
Z\left( t,t;z,v\right) =z\;\;,\;\;V\left( t,t;z,v\right) =v   \label{W4}
\end{equation}
with $s>t.$

Using Taylor's expansion we can write: 
\begin{align*}
\frac{dZ}{ds} & =-sE\left( z+vs,s\right) \\
& +O\left( s\left\| E\right\| e^{-\gamma s}\left| Z-z\right| +s^{2}\left\|
E\right\| e^{-\gamma s}\left| V-v\right| +\left\| E\right\| se^{-\frac{\gamma%
}{2}s}\exp\left( -\frac{a}{\left\| E\right\| }e^{\gamma t}\right) \right) \\
\frac{dV}{ds} & =E\left( z+vs,s\right) \\
& +O\left( \left\| E\right\| e^{-\gamma s}\left| Z-z\right| +s\left\|
E\right\| e^{-\gamma s}\left| V-v\right| +\left\| E\right\| e^{-\frac{\gamma%
}{2}s}\exp\left( -\frac{a}{\left\| E\right\| }e^{\gamma t}\right) \right)
\end{align*}
for some $a>0.$ The first two terms of the remainder arise in making a
Taylor expansion for the range of values where $s\left| V-v\right| \leq1.$
The last term is the contribution from the region where $s\left| V-v\right|
>1.$ Since $\left| V-v\right| \leq CM\left\| E\right\| e^{-\gamma t}$ this
requires $s$ huge, and due to the exponential decay of $E$ the estimate
follows. Taking into account the estimates for the differences $\left|
Z-z\right| ,\;\left| V-v\right| $ that are similar to the estimates obtained
for $Z_{\infty},\;V_{\infty}$ it then follows: 
\begin{align}
Z_{\infty}\left( t;z,v\right) -z & =-\int_{t}^{\infty}sE\left( z+vs,s\right)
ds+O\left( \left( t+1\right) ^{2}\left\| E\right\| ^{2}e^{-2\gamma t}\right)
\label{W4a} \\
V_{\infty}\left( t;z,v\right) -v & =\int_{t}^{\infty}E\left( z+vs,s\right)
ds+O\left( \left( t+1\right) \left\| E\right\| ^{2}e^{-2\gamma t}\right) 
\label{W4b}
\end{align}

We can then write the following equation for $E\left( z,t\right) :$%
\begin{align}
E_{z}\left( z,t\right) & =\int_{t}^{\infty}ds\int_{-\infty}^{\infty
}f_{e,v}\left( w\right) E\left( z+ws-wt,s\right)
dw+\int_{-\infty}^{\infty}g_{\infty}\left( z-wt,w\right) dw  \notag \\
& -\int_{t}^{\infty}sds\int_{-\infty}^{\infty}\frac{\partial g_{\infty}}{%
\partial z}\left( z-wt,w\right) E\left( z+ws-wt,s\right) dw+  \notag \\
& +\int_{t}^{\infty}ds\int_{-\infty}^{\infty}\frac{\partial g_{\infty}}{%
\partial w}\left( z-wt,w\right) E\left( z+ws-wt,s\right) dw+\tilde {R}\left(
z,t\right)   \label{W7a}
\end{align}
where\ we remark that the definition of $\tilde{R}$ is: 
\begin{align}
\tilde{R}\left( z,t\right) & =\int_{-\infty}^{\infty}\left[ f_{e}\left(
V_{\infty}\left( t;z-wt,w\right) \right) -f_{e}\left( w\right) \right] dw+ 
\notag \\
& +\int_{-\infty}^{\infty}g_{\infty}\left( Z_{\infty}\left( t;z-wt,w\right)
,V_{\infty}\left( t;z-wt,w\right) \right) dw-  \notag \\
& -\int_{t}^{\infty}ds\int_{-\infty}^{\infty}f_{e,v}\left( w\right) E\left(
z+ws-wt,s\right) dw-\int_{-\infty}^{\infty}g_{\infty}\left( z-wt,w\right) dw+
\notag \\
& +\int_{t}^{\infty}sds\int_{-\infty}^{\infty}\frac{\partial g_{\infty}}{%
\partial z}\left( z-wt,w\right) E\left( z+ws-wt,s\right) dw-  \notag \\
& -\int_{t}^{\infty}ds\int_{-\infty}^{\infty}\frac{\partial g_{\infty}}{%
\partial w}\left( z-wt,w\right) E\left( z+ws-wt,s\right) dw   \label{W7}
\end{align}

Using Assumptions (A), expanding the arguments of $f_{e},\;g_{\infty}$ using
Taylor's expansion as well as the estimates (\ref{W4a}), (\ref{W4b}) we
obtain: 
\begin{equation}
\left| \tilde{R}\left( z,t\right) \right| \leq C\left\| E\right\| ^{2}\left(
t+1\right) ^{2}e^{-2\gamma t}   \label{W7new}
\end{equation}

We now transform the third term on the right side as suggested above, using
integrations by parts. It then follows that (see (\ref{N9})): 
\begin{align*}
E_{z}\left( z,t\right) & =\int_{t}^{\infty}ds\int_{-\infty}^{\infty
}f_{e,v}\left( w\right) E\left( z-wt+ws,s\right)
dw+\int_{-\infty}^{\infty}g_{\infty}\left( z-wt,w\right) dw \\
& +\int_{t}^{\infty}\left( t-s\right) \frac{s}{t}ds\int_{-\infty}^{\infty
}g_{\infty}\left( z-wt,w\right) E_{z}\left( z-wt+ws,s\right) dw+ \\
& +\int_{t}^{\infty}ds\int_{-\infty}^{\infty}\left( 1-\frac{s}{t}\right) 
\frac{\partial g_{\infty}}{\partial w}\left( z-wt,w\right) E\left(
z-wt+ws,s\right) dw+\tilde{R}\left( z,t\right)
\end{align*}

A crucial estimate is the following: 
\begin{align}
& \left| \int_{t}^{\infty}\left( t-s\right) \frac{s}{t}ds\int_{-\infty
}^{\infty}g_{\infty}\left( z-wt,w\right) E_{z}\left( z-wt+ws,s\right)
dw\right|  \label{W7new2} \\
& \leq C\varepsilon\left\| E\right\| e^{-\gamma t}\int_{t}^{\infty}\left(
s-t\right) \frac{s}{t}e^{-\gamma\left( s-t\right) }ds\leq C\varepsilon
\left\| E\right\| e^{-\gamma t}  \notag
\end{align}
where $\varepsilon$ comes from $g_{\infty}.$ A similar estimate can be
obtained for the fourth term. The problem is then ready for a fixed point
argument in the form 
\begin{equation}
E_{z}\left( z,t\right)
-\int_{t}^{\infty}ds\int_{-\infty}^{\infty}f_{e,v}\left( w\right) E\left(
z-wt+ws,s\right) dw=\psi\left( z,t\right) +L\left( z,t\right) +\tilde{R}%
\left( z,t\right)   \label{W5}
\end{equation}
where: 
\begin{align}
\psi\left( z,t\right) & =\int_{-\infty}^{\infty}g_{\infty}\left(
z-wt,w\right) dw  \label{W6a} \\
L\left( z,t\right) & =\int_{t}^{\infty}\left( t-s\right) \frac{s}{t}%
ds\int_{-\infty}^{\infty}g_{\infty}\left( z-wt,w\right) E_{z}\left(
z-wt+ws,s\right) dw+  \label{W6b} \\
& \int_{t}^{\infty}ds\int_{-\infty}^{\infty}\left( 1-\frac{s}{t}\right) 
\frac{\partial g_{\infty}}{\partial w}\left( z-wt,w\right) E\left(
z-wt+ws,s\right) dw  \notag
\end{align}

More precisely, the fixed point scheme is the following one. Given $E$ with $%
\left\| E\right\| <\infty$ we define $L$ as in (\ref{W6b}) and $\tilde{R}$
as in (\ref{W7}). We then define $\mathcal{T}\left( E\right) $ as the
solution of (\ref{W5}) obtained by means of Proposition \ref{regestimate}.

In particular, using Lemma \ref{expdecay}, (\ref{W7new}), (\ref{W7new2}) and
Proposition \ref{regestimate} we have the following estimate: 
\begin{equation}
\left\| \mathcal{T}\left( E\right) \right\| \leq
C_{g_{\infty}}\varepsilon+C\varepsilon\left\| E\right\| +C\left\| E\right\|
^{2}   \label{W1}
\end{equation}
where $C_{g_{\infty}}$ is a constant of order one that comes from $g_{\infty}
$ and depends only on $B,\;\alpha.$

A consequence of (\ref{W1}) is that the ball $\left\| E\right\|
\leq2C_{g_{\infty}}\varepsilon$ is transformed by means of the operator $%
\mathcal{T}$ in a set contained in a ball with radius: 
\begin{equation*}
C_{g_{\infty}}\varepsilon+2CC_{g_{\infty}}\varepsilon^{2}+4C\left(
C_{g_{\infty}}\right) ^{2}\varepsilon^{2}
\end{equation*}
and this number is strictly less than $2C_{g_{\infty}}\varepsilon$ if $%
\varepsilon$ is sufficiently small.

Finally we will prove that the operator $\mathcal{T}$ is contractive in the
ball $\left\| E\right\| \leq2C_{g_{\infty}}\varepsilon$ if $\varepsilon$ is
sufficiently small. To this end we denote as $Z_{1}\left( s,t;z,v\right)
,\;V_{1}\left( s,t;z,v\right) ,\allowbreak\;Z_{2}\left( s,t;z,v\right)
,\;V_{2}\left( s,t;z,v\right) $ the evolution of the characteristics defined
by means of (\ref{W2})-(\ref{W4}) with electric fields $E_{1},\;E_{2}$
respectively with $\left\| E_{1}\right\| \leq2C_{g_{\infty}}\varepsilon
,\;\left\| E_{2}\right\| \leq2C_{g_{\infty}}\varepsilon.$ We will use also
the notation $Z_{1,\infty}\left( t;z,v\right) ,$ $V_{1,\infty}\left(
t;z,v\right) ,$ $Z_{2,\infty}\left( t;z,v\right) ,$ $V_{2,\infty}\left(
t;z,v\right) ,$ $L_{1}\left( z,t\right) ,$ $L_{2}\left( z,t\right) ,$ $%
\tilde{R}_{1}\left( z,t\right) ,$ $\tilde{R}_{2}\left( z,t\right) $ to
denote the corresponding functions associated with $E_{1},\;E_{2}.$

Notice that (\ref{W2})-(\ref{W4}) imply: 
\begin{align*}
\frac{\partial\left( Z_{1}-Z_{2}\right) }{\partial s} & =-s\left[
E_{1}\left( Z_{1}+V_{1}s,s\right) -E_{2}\left( Z_{2}+V_{2}s,s\right) \right]
\\
\frac{\partial\left( V_{1}-V_{2}\right) }{\partial s} & =\left[ E_{1}\left(
Z_{1}+V_{1}s,s\right) -E_{2}\left( Z_{2}+V_{2}s,s\right) \right]
\end{align*}
with: 
\begin{align}
\left( Z_{1}-Z_{2}\right) \left( t,t;z,v\right) & =0  \label{Z1} \\
\left( V_{1}-V_{2}\right) \left( t,t;z,v\right) & =0  \notag
\end{align}

We can write the equations as: 
\begin{align}
\frac{\partial\left( Z_{1}-Z_{2}\right) }{\partial s} & =-s\left[ a\left(
s\right) \left( Z_{1}-Z_{2}\right) +b\left( s\right) \left(
V_{1}-V_{2}\right) +\left( E_{1}\left( Z_{2}+V_{2}s,s\right) -E_{2}\left(
Z_{2}+V_{2}s,s\right) \right) \right]  \label{Z2} \\
\frac{\partial\left( V_{1}-V_{2}\right) }{\partial s} & =\left[ a\left(
s\right) \left( Z_{1}-Z_{2}\right) +b\left( s\right) \left(
V_{1}-V_{2}\right) +\left( E_{1}\left( Z_{2}+V_{2}s,s\right) -E_{2}\left(
Z_{2}+V_{2}s,s\right) \right) \right]   \label{Z3}
\end{align}
where: 
\begin{align*}
a\left( s\right) & =\frac{E_{1}\left( Z_{1}+V_{1}s,s\right) -E_{1}\left(
Z_{2}+V_{1}s,s\right) }{Z_{1}-Z_{2}} \\
b\left( s\right) & =\frac{E_{1}\left( Z_{2}+V_{1}s,s\right) -E_{1}\left(
Z_{2}+V_{2}s,s\right) }{V_{1}-V_{2}}
\end{align*}

Notice that: 
\begin{align}
\left| a\left( s\right) \right| +\left| b\left( s\right) \right| &
\leq\left\| E_{1}\right\| e^{-\gamma s}  \label{Z4} \\
\left| E_{1}\left( Z_{2}+V_{2}s,s\right) -E_{2}\left( Z_{2}+V_{2}s,s\right)
\right| & \leq\left\| E_{1}-E_{2}\right\| e^{-\gamma s}   \label{Z5}
\end{align}

Then, combining (\ref{Z1})-(\ref{Z5}) we obtain, using a Gronwall type
argument: 
\begin{align*}
\left| \left( Z_{1}-Z_{2}\right) \left( s,t;z,v\right) \right| & \leq
C\left\| E_{1}-E_{2}\right\| \left( t+1\right) e^{-\gamma t} \\
\left| \left( V_{1}-V_{2}\right) \left( s,t;z,v\right) \right| & \leq
C\left\| E_{1}-E_{2}\right\| e^{-\gamma t}
\end{align*}
whence: 
\begin{align}
\left| \left( Z_{1,\infty}-Z_{2,\infty}\right) \left( t;z,v\right) \right| &
\leq C\left\| E_{1}-E_{2}\right\| \left( t+1\right) e^{-\gamma t}\;\;,
\label{Z6} \\
\;\;\left| \left( V_{1,\infty}-V_{2,\infty}\right) \left( t;z,v\right)
\right| & \leq C\left\| E_{1}-E_{2}\right\| e^{-\gamma t}  \notag
\end{align}

Let us denote as $\tilde{E}_{1},\;\tilde{E}_{2}$ respectively the solutions
of (\ref{W5}) with sources $L_{1},\;\tilde{R}_{1}$ and $L_{2},\;\tilde{R}_{2}
$ respectively.

Notice that: 
\begin{align}
& \left( \tilde{E}_{1}-\tilde{E}_{2}\right) _{z}\left( z,t\right)
-\int_{t}^{\infty}ds\int_{-\infty}^{\infty}f_{e,v}\left( w\right) \left( 
\tilde{E}_{1}-\tilde{E}_{2}\right) \left( z-wt+ws,s\right) dw  \label{Z6a} \\
& =\left( L_{1}-L_{2}\right) \left( z,t\right) +\left( \tilde{R}_{1}-\tilde{R%
}_{2}\right) \left( z,t\right)  \notag
\end{align}

Notice also that, arguing as in the estimates of $L$ above (cf. (\ref{W6b})
and (\ref{W7new2})): 
\begin{equation}
\left| \left( L_{1}-L_{2}\right) \left( z,t\right) \right| \leq
C\varepsilon\left\| E_{1}-E_{2}\right\| e^{-\gamma t}\;\;,\;\;z\in \mathbb{%
R\;\;},\;\;t\geq0   \label{Z7}
\end{equation}

We also need to estimate the difference $\left( \tilde{R}_{1}-\tilde{R}%
_{2}\right) .$ Therefore we need to estimate: 
\begin{align}
\left( \tilde{R}_{1}-\tilde{R}_{2}\right) \left( z,t\right) &
=\int_{-\infty}^{\infty}\left[ f_{e}\left( V_{1,\infty}\left(
t;z-wt,w\right) \right) -f_{e}\left( V_{2,\infty}\left( t;z-wt,w\right)
\right) \right] dw+  \notag \\
& +\int_{-\infty}^{\infty}\left[ g_{\infty}\left( Z_{1,\infty}\left(
t;z-wt,w\right) ,V_{1,\infty}\left( t;z-wt,w\right) \right) -\right.  \notag
\\
& \left. -g_{\infty}\left( Z_{2,\infty}\left( t;z-wt,w\right)
,V_{2,\infty}\left( t;z-wt,w\right) \right) \right] dw+  \notag \\
& -\int_{t}^{\infty}ds\int_{-\infty}^{\infty}f_{e,v}\left( w\right) \left(
E_{1}\left( z+ws-wt,s\right) -E_{2}\left( z+ws-wt,s\right) \right) dw+ 
\notag \\
& +\int_{t}^{\infty}sds\int_{-\infty}^{\infty}\frac{\partial g_{\infty}}{%
\partial z}\left( z-wt,w\right) \left( E_{1}\left( z+ws-wt,s\right)
-E_{2}\left( z+ws-wt,s\right) \right) dw-  \notag \\
& -\int_{t}^{\infty}ds\int_{-\infty}^{\infty}\frac{\partial g_{\infty}}{%
\partial w}\left( z-wt,w\right) \left( E_{1}\left( z+ws-wt,s\right)
-E_{2}\left( z+ws-wt,s\right) \right) dw   \label{Z8}
\end{align}

Taking into account (\ref{A2}) and (\ref{Z6}) we can estimate the second,
fourth and fifth terms in (\ref{Z8}) as $C\varepsilon\left\|
E_{1}-E_{2}\right\| e^{-\gamma t}.$ Then: 
\begin{align}
\left| \left( \tilde{R}_{1}-\tilde{R}_{2}\right) \left( z,t\right) \right| &
\leq\left| \int_{-\infty}^{\infty}\left[ f_{e}\left( V_{1,\infty}\left(
t;z-wt,w\right) \right) -f_{e}\left( V_{2,\infty }\left( t;z-wt,w\right)
\right) \right] dw-\right.  \notag \\
& \left. -\int_{t}^{\infty}ds\int_{-\infty}^{\infty}f_{e,v}\left( w\right)
\left( E_{1}\left( z+ws-wt,s\right) -E_{2}\left( z+ws-wt,s\right) \right)
dw\right| +  \notag \\
& +C\varepsilon\left\| E_{1}-E_{2}\right\| e^{-\gamma t}   \label{Z9}
\end{align}

On the other hand (\ref{W3}), (\ref{W4}) imply: 
\begin{equation}
V_{i,\infty}\left( t;z-wt,w\right) -w=\int_{t}^{\infty}E_{i}\left(
Z_{i}\left( s\right) +sV_{i}\left( s\right) ,s\right) ds\;\;,\;\;i=1,2 
\label{Z10}
\end{equation}
where, for simplicity $Z_{i}\left( s\right) =Z_{i}\left( s,t;z-wt,w\right)
,\;V_{i}\left( s\right) =V_{i}\left( s,t;z-wt,w\right) \;,\;\;i=1,2.$

Using Taylor's Theorem and (\ref{Z6}) we can estimate the difference\ 
\newline
\ $\left[ f_{e}\left( V_{1,\infty}\left( t;z-wt,w\right) \right)
-f_{e}\left( V_{2,\infty}\left( t;z-wt,w\right) \right) \right] $ as: 
\begin{align*}
& \left[ f_{e}\left( V_{1,\infty}\left( t;z-wt,w\right) \right) -f_{e}\left(
V_{2,\infty}\left( t;z-wt,w\right) \right) \right] \\
& =f_{e,v}\left( V_{2,\infty}\left( t;z-wt,w\right) \right) \left(
V_{1,\infty}\left( t;z-wt,w\right) -V_{2,\infty}\left( t;z-wt,w\right)
\right) + \\
& +O\left( \left\| E_{1}-E_{2}\right\| ^{2}e^{-2\gamma t}\right)
\end{align*}

Using this estimate as well as (\ref{Z6}), (\ref{Z10}) and the fact that $%
\left\| E_{i}\right\| \leq2C_{g_{\infty}}\varepsilon,\;i=1,2$ we obtain: 
\begin{align*}
& \left[ f_{e}\left( V_{1,\infty}\left( t;z-wt,w\right) \right) -f_{e}\left(
V_{2,\infty}\left( t;z-wt,w\right) \right) \right] \\
& =f_{e,v}\left( w\right) \left( V_{1,\infty}\left( t;z-wt,w\right)
-V_{2,\infty}\left( t;z-wt,w\right) \right) + \\
& +O\left( \left\| E_{1}-E_{2}\right\| ^{2}e^{-2\gamma t}\right) +O\left(
\varepsilon\left\| E_{1}-E_{2}\right\| e^{-\gamma t}\right)
\end{align*}

Taking into account (\ref{Z10}) it then follows that: 
\begin{align*}
& \left| \int_{-\infty}^{\infty}\left[ f_{e}\left( V_{1,\infty}\left(
t;z-wt,w\right) \right) -f_{e}\left( V_{2,\infty}\left( t;z-wt,w\right)
\right) \right] dw-\right. \\
& \left. -\int_{t}^{\infty}ds\int_{-\infty}^{\infty}f_{e,v}\left( w\right)
\left( E_{1}\left( z+ws-wt,s\right) -E_{2}\left( z+ws-wt,s\right) \right)
dw\right| \\
& \leq\left| \int_{t}^{\infty}ds\int_{-\infty}^{\infty}f_{e,v}\left(
w\right) \left( \left[ E_{1}\left( Z_{1}\left( s\right) +sV_{1}\left(
s\right) ,s\right) -E_{1}\left( z+ws-wt,s\right) \right] -\right. \right. \\
& \left. -\left. \left[ E_{2}\left( Z_{2}\left( s\right) +sV_{2}\left(
s\right) ,s\right) -E_{2}\left( z+ws-wt,s\right) \right] \right) \right| + \\
& +C\left( \left\| E_{1}-E_{2}\right\| ^{2}e^{-2\gamma t}+\varepsilon
\left\| E_{1}-E_{2}\right\| e^{-\gamma t}\right)
\end{align*}

Notice that 
\begin{align*}
& \left( E_{1}\left( Z_{1}\left( s\right) +sV_{1}\left( s\right) ,s\right)
-E_{1}\left( z+ws-wt,s\right) \right) \\
& -\left( E_{2}\left( Z_{2}\left( s\right) +sV_{2}\left( s\right) ,s\right)
-E_{2}\left( z+ws-wt,s\right) \right) \\
& =\left( E_{1}\left( Z_{1}\left( s\right) +sV_{1}\left( s\right) ,s\right)
-E_{1}\left( Z_{2}\left( s\right) +sV_{2}\left( s\right) ,s\right) \right) +
\\
& +\left( E_{1}\left( Z_{2}\left( s\right) +sV_{2}\left( s\right) ,s\right)
-E_{1}\left( z+ws-wt,s\right) \right) \\
& -\left( E_{2}\left( Z_{2}\left( s\right) +sV_{2}\left( s\right) ,s\right)
-E_{2}\left( z+ws-wt,s\right) \right)
\end{align*}

The term $\left( E_{1}\left( Z_{1}\left( s\right) +sV_{1}\left( s\right)
,s\right) -E_{1}\left( Z_{2}\left( s\right) +sV_{2}\left( s\right) ,s\right)
\right) $ can be estimated as: 
\begin{equation*}
\left| \int_{Z_{2}\left( s\right) +sV_{2}\left( s\right) }^{Z_{1}\left(
s\right) +sV_{1}\left( s\right) }E_{1,z}\left( \xi,s\right) d\xi\right| \leq
C\left\| E_{1}\right\| e^{-\gamma s}\left[ \left| Z_{1}\left( s\right)
-Z_{2}\left( s\right) \right| +s\left| V_{1}\left( s\right) -V_{2}\left(
s\right) \right| \right] 
\end{equation*}
whence, using $\left\| E_{i}\right\| \leq2C_{g_{\infty}}\varepsilon$: 
\begin{equation}
\left| E_{1}\left( Z_{1}\left( s\right) +sV_{1}\left( s\right) ,s\right)
-E_{1}\left( Z_{2}\left( s\right) +sV_{2}\left( s\right) ,s\right) \right|
\leq C\varepsilon e^{-\gamma s}e^{-\gamma t}\left( 1+s\right) \left( \left\|
E_{1}-E_{2}\right\| \right)   \label{Z11}
\end{equation}

On the other hand 
\begin{align*}
& \left( E_{1}\left( Z_{2}\left( s\right) +sV_{2}\left( s\right) ,s\right)
-E_{1}\left( z+ws-wt,s\right) \right) - \\
& -\left( E_{2}\left( Z_{2}\left( s\right) +sV_{2}\left( s\right) ,s\right)
-E_{2}\left( z+ws-wt,s\right) \right) \\
& =\int_{z+ws-wt}^{Z_{2}\left( s\right) +sV_{2}\left( s\right) }\left[
E_{1,z}\left( \xi,s\right) -E_{2,z}\left( \xi,s\right) \right] d\xi
\end{align*}
thus we have 
\begin{align}
& \left| \left( E_{1}\left( Z_{2}\left( s\right) +sV_{2}\left( s\right)
,s\right) -E_{1}\left( z+ws-wt,s\right) \right) -\right.  \notag \\
& \left. -\left( E_{2}\left( Z_{2}\left( s\right) +sV_{2}\left( s\right)
,s\right) -E_{2}\left( z+ws-wt,s\right) \right) \right|  \notag \\
& \leq\left\| E_{1}-E_{2}\right\| e^{-\gamma s}\left[ \left| Z_{2}\left(
s\right) -\left( z+ws-wt\right) \right| +s\left| V_{2}\left( s\right)
-w\right| \right]  \notag \\
& \leq Ce^{-\gamma s}e^{-\gamma t}\left( 1+s\right) \left\|
E_{1}-E_{2}\right\| \left\| E_{2}\right\|   \label{Z12}
\end{align}

Then: 
\begin{align*}
& \left| \int_{-\infty}^{\infty}\left[ f_{e}\left( V_{1,\infty}\left(
t;z-wt,w\right) \right) -f_{e}\left( V_{2,\infty}\left( t;z-wt,w\right)
\right) \right] dw-\right. \\
& \left. -\int_{t}^{\infty}ds\int_{-\infty}^{\infty}f_{e,v}\left( w\right)
\left( E_{1}\left( z+ws-wt,s\right) -E_{2}\left( z+ws-wt,s\right) \right)
dw\right| \\
& \leq C\varepsilon\left\| E_{1}-E_{2}\right\| e^{-\gamma t}
\end{align*}
and therefore: 
\begin{equation}
\left| \left( \tilde{R}_{1}-\tilde{R}_{2}\right) \left( z,t\right) \right|
\leq C\varepsilon\left\| E_{1}-E_{2}\right\| e^{-\gamma t}   \label{Z12a}
\end{equation}

Combining (\ref{Z6}), (\ref{Z6a}), (\ref{Z12a}) and Proposition \ref
{regestimate} we obtain: 
\begin{equation*}
\left\| \mathcal{T}\left( E_{1}\right) -\mathcal{T}\left( E_{2}\right)
\right\| \leq C\varepsilon\left\| E_{1}-E_{2}\right\| 
\end{equation*}
for $\left\| E_{i}\right\| \leq2C_{g_{\infty}}\varepsilon,\;i=1,2,$ thus the
desired contractibility of $\mathcal{T}$ follows. This concludes the proof.
\end{proof}

\bigskip

\textbf{Acknowledgements:} This work has been partially supported by the
Korean Research Fund KRF 2008-314-C00023 and Grant of the Spanish Ministery
of Science and Innovation DGES MTM2007-61755. JJLV also thanks the
hospitality of Universidad Complutense where part of this research has been
made.


\bigskip

\begin{thebibliography}{99}
\bibitem{Ahlfors}  L. Ahlfors, Complex Analysis. Mc. Graw Hill, third
edition, (1979).

\bibitem{Arnold}  V. I. Arnold, Mathematical Methods of Classical Mechanics.
Springer Verlag, second edition (1997).

\bibitem{BohmGross}  D. Bohm, E. P. Gross, Theory of plasma oscillations. A.
Origin of medium-like behaviour. Phys. Rev. 75, 12, 1851-1864, (1949).

\bibitem{CM}  E. Caglioti, C. Maffei, Time asymptotics for solutions of
Vlasov-Poisson equation in a circle. J. Stat. Phys. 92, n 1/2, (1998).

\bibitem{Degond}  P. Degond, Spectral Theory of the linearized
Vlasov-Poisson equation. Trans. AMS, vol. 294, 2, 435-453, (1986).

\bibitem{Glassey}  R. Glassey, The Cauchy problem in kinetic theory, SIAM:
Philadelphia, PA, 1996.

\bibitem{GuoStrauss}  Y. Guo, W. Strauss, Instability of periodic BGK
equilibria. Comm. Pure Appl. Math. 48, no. 8, 861--894, (1995).

\bibitem{Jackson}  J. D. Jackson, Longitudinal plasma oscillations. J. Nucl.
Energy, Part C: Plasma Physics, vol. 1, 171-189, (1960).

\bibitem{Landau}  L. Landau, On the vibration of the electronic plasma.
Journal of Physics, v. X, n. 1, 25-34, (1946).

\bibitem{LP}  P. L. Lions, B. Perthame, Propagation of moments and
regularity of solutions for the 3-dimensional Vlasov-Poisson system. Invent.
Math. \textbf{105}, 415-430, (1991).

\bibitem{Pfaffelmoser}  K. Pfaffelmoser, Global classical solutions of the
Vlasov-Poisson system in three dimensions for general initial data. J.
Differential Equations \textbf{95}, 281-303, (1992).

\bibitem{Saenz}  A. W. Saenz, Long-time behaviour of the electric potential
and stability in the linearized Vlasov theory. Journal of Mathematical
Physics, vol. 6, n. 6, 859-875, (1965).

\bibitem{Schaefer}  R. Glassey, J. Schaeffer, Time decay for solutions to
the linearized Vlasov equation. Transport Theory Statist. Phys. 23, no. 4,
411--453, (1994).

\bibitem{Schaffer2}  J. Schaeffer, Global existence of smooth solutions to
the Vlasov-Poisson system in three dimensions, Comm. Partial Differential
Equations, 16, 1313--1335, (1991).

\bibitem{Vlasov}  A. Vlasov, On the kinetic theory of an assembly of
particles with collective interaction, J. Phys. (USSR), (1945).

\bibitem{Guo}  T. Zhou, Y. Guo, C.-W. Shu, Numerical study on Landau
damping, Physica D, 157, 322-333, (2001).
\end{thebibliography}
\end{document}